\documentclass[11pt,a4paper]{amsart}
\usepackage{amsthm,amsmath,amsfonts,graphicx}
\usepackage[sort&compress,numbers]{natbib}
\usepackage[lmargin=32mm,rmargin=32mm,tmargin=32mm,bmargin=32mm]{geometry}

\newcommand{\spacing}[1]{
\renewcommand{\baselinestretch}{#1}
\setlength{\footnotesep}{\baselinestretch\footnotesep}}
\spacing{1.15}

\theoremstyle{plain}
\newtheorem{thm}{Theorem}[section]
\newtheorem{lem}[thm]{Lemma}
\newtheorem{cor}[thm]{Corollary}
\newtheorem{conj}[thm]{Conjecture}
\newtheorem{prop}[thm]{Proposition}
\newtheorem{claim}[thm]{Claim}

\newcommand{\Figure}[4][htb]{
\begin{figure}[#1]
    \vspace*{1ex}
    \begin{center}#3\end{center}
    \vspace*{-1ex}
    \caption{#4}
    \label{fig:#2}
\end{figure}}

\newcommand{\isomorphic}{\protect\cong}

\newcommand{\COR}[2]{\begin{cor}\label{cor:#1}#2\end{cor}}
\newcommand{\THM}[2]{\begin{thm}\label{thm:#1}#2\end{thm}}

\newcommand{\LEM}[2]{\begin{lem}\label{lem:#1}#2\end{lem}}

\newcommand{\PROPOS}[3]{\begin{prop}[#1]\label{prop:#2}#3\end{prop}}
\newcommand{\LEMMA}[3]{\begin{lem}[#1]\label{lem:#2}#3\end{lem}}
\newcommand{\SECT}[2]{\section{#2}\label{sec:#1}}
\newcommand{\SUBSECT}[2]{\subsection{#2}\label{sec:#1}}
\newcommand{\EQN}[2]{\begin{equation}\label{eqn:#1}#2\end{equation}}
\newcommand{\eqn}[1]{\begin{equation*}#1\end{equation*}}
\newcommand{\PROOF}[1]{\begin{proof}#1\end{proof}}

\newcommand{\corref}[1]{Corollary~\ref{cor:#1}}
\newcommand{\thmref}[1]{Theorem~\ref{thm:#1}}
\newcommand{\lemref}[1]{Lemma~\ref{lem:#1}}
\newcommand{\propref}[1]{Proposition~\ref{prop:#1}}

\newcommand{\twolemref}[2]{Lemmas~\ref{lem:#1} and \ref{lem:#2}}
\newcommand{\threelemref}[3]{Lemmas~\ref{lem:#1}, \ref{lem:#2} and \ref{lem:#3}}

\newcommand{\secref}[1]{Section~\ref{sec:#1}}

\newcommand{\eqnref}[1]{Equation~\eqref{eqn:#1}}

\newcommand{\figref}[1]{Figure~\ref{fig:#1}}
\newcommand{\twofigref}[2]{Figures~\ref{fig:#1} and \ref{fig:#2}}
\newcommand{\manyfigref}[2]{Figures~\ref{fig:#1}--\ref{fig:#2}}

\DeclareMathOperator{\diam}{diam}
\DeclareMathOperator{\dist}{dist}
\DeclareMathOperator{\ecc}{ecc}

\newcommand{\bracket}[1]{\ensuremath{\protect\left(#1\right)}}
\newcommand{\BRACKET}[1]{\ensuremath{\protect\left(#1\right)}}

\newcommand{\FLOOR}[1]{\ensuremath{\protect\left\lfloor#1\right\rfloor}}
\newcommand{\floor}[1]{\ensuremath{\protect\lfloor#1\rfloor}}
\newcommand{\ceil}[1]{\ensuremath{\protect\lceil#1\rceil}}
\newcommand{\CEIL}[1]{\ensuremath{\protect\left\lceil#1\right\rceil}}

\newcommand{\G}[1][\beta,D]{\ensuremath{\mathcal{G}_{#1}}}

\begin{document}

\title[Metric Dimension and Diameter]{Extremal Graph Theory for\\ Metric Dimension and Diameter}

\thanks{The research of Carmen Hernando, Merc{\`e} Mora, Carlos Seara, and David Wood is supported by the projects MEC MTM2006-01267 and DURSI 2005SGR00692.  
The research of Ignacio Pelayo is supported by the projects MTM2005-08990-C02-01 and SGR2005-00412. The research of David Wood is supported by a Marie Curie Fellowship of the European Community under contract 023865.}

\author[Hernando]{Carmen~Hernando}
\address{Departament de Matem{\`a}tica Aplicada I, Universitat Polit{\`e}cnica de Catalunya, Barcelona, Spain}
\email{carmen.hernando@upc.edu}

\author[Mora]{Merc{\`e}~Mora}
\address{Departament de Matem{\`a}tica Aplicada II, Universitat Polit{\`e}cnica de Catalunya, Barcelona, Spain}
\email{merce.mora@upc.edu}

\author[Pelayo]{Ignacio~M.~Pelayo}
\address{Departament de Matem{\`a}tica Aplicada III, Universitat Polit{\`e}cnica de Catalunya, Barcelona, Spain}
\email{ignacio.m.pelayo@upc.edu}

\author[Seara]{Carlos~Seara}
\address{Departament de Matem{\`a}tica Aplicada II, Universitat Polit{\`e}cnica de Catalunya, Barcelona, Spain}
\email{carlos.seara@upc.edu}

\author[Wood]{David~R.~Wood}
\address{Departament de Matem{\`a}tica Aplicada II, Universitat Polit{\`e}cnica de Catalunya, Barcelona, Spain}
\email{david.wood@upc.edu}

\keywords{graph, distance, resolving set, metric dimension, metric basis, diameter, order}

\subjclass[2000]{05C12 (distance in graphs), 05C35 (extremal graph theory)}

\begin{abstract}
A set of vertices $S$ \emph{resolves} a connected graph $G$ if every vertex is uniquely determined by its vector of distances to the vertices in $S$. The \emph{metric dimension} of $G$ is the minimum cardinality of a resolving set of $G$. Let $\mathcal{G}_{\beta,D}$ be the set of graphs with metric dimension $\beta$ and diameter $D$. It is well-known that the minimum order of a graph in $\mathcal{G}_{\beta,D}$ is exactly $\beta+D$. The first contribution of this paper is to characterise the graphs in $\mathcal{G}_{\beta,D}$ with order $\beta+D$ for all values of $\beta$ and $D$. Such a characterisation was previously only known for $D\leq2$ or $\beta\leq1$. The second contribution is to determine the maximum order of a graph in $\mathcal{G}_{\beta,D}$ for all values of $D$ and $\beta$. Only a weak upper bound was previously known.
\end{abstract}

\maketitle

\SECT{Introduction}{Introduction}

Let $G$ be a connected graph\footnote{Graphs in this paper are finite, undirected, and simple. The vertex set and edge set of a graph $G$ are denoted by $V(G)$ and $E(G)$. For vertices $v,w\in V(G)$, we write $v\sim w$ if $vw\in E(G)$, and $v\not\sim w$ if $vw\not\in E(G)$. For $S\subseteq V(G)$, let $G[S]$ be the subgraph of $G$ induced by $S$. That is, $V(G[S])=S$ and $E(G[S])=\{vw\in E(G):v\in S,w\in S\}$). For $S\subseteq V(G)$, let $G\setminus S$ be the graph $G[V(G)\setminus S]$. For $v\in V(G)$, let $G\setminus v$ be the graph $G\setminus\{v\}$. Suppose that $G$ is connected. The \emph{distance} between vertices $v,w\in V(G)$, denoted by $\dist_G(v,w)$, is the length (that is, the number of edges) in a shortest path between $v$ and $w$ in $G$. The \emph{eccentricity} of a vertex $v$ in $G$ is $\ecc_G(v):=\max\{\dist_G(v,w):w\in V(G)\}$. We drop the subscript $G$ from these notations if the graph $G$ is clear from the context. The \emph{diameter} of $G$ is $\diam(G):=\max\{\dist(v,w):v,w\in V(G)\}=\max\{\ecc(v):v\in V(G)\}$. For integers $a\leq b$, let $[a,b]:=\{a,a+1,\dots,b\}$.}. A vertex $x\in V(G)$ \emph{resolves}\footnote{It will be convenient to also use the following definitions for a connected graph $G$. A vertex $x\in V(G)$ \emph{resolves} a set of vertices $T\subseteq V(G)$ if $x$ resolves every pair of distinct vertices in $T$. A set of vertices $S\subseteq V(G)$ \emph{resolves} a set of vertices $T\subseteq V(G)$ if for every pair of distinct vertices $v,w\in T$, there exists a vertex $x\in S$ that resolves $v,w$.} a pair of vertices $v,w\in V(G)$ if $\dist(v,x)\ne \dist(w,x)$. A set of vertices $S\subseteq V(G)$ \emph{resolves} $G$, and $S$ is a \emph{resolving set} of $G$, if every pair of distinct vertices of $G$ are resolved by some vertex in $S$. Informally, $S$ resolves $G$ if every vertex of $G$ is uniquely determined by its vector of distances to the vertices in $S$. A resolving set $S$ of $G$ with the minimum cardinality is a \emph{metric basis} of $G$, and $|S|$ is the \emph{metric dimension} of $G$, denoted by $\beta(G)$.

Resolving sets in general graphs were first defined by \citet{Slater75} and \citet{HM-AC76}. Resolving sets have since been widely investigated \citep{SSH02, SZ-CN04, SZ-IJMMS04, CZ-CN03, Caceres-etal, PO-UM06, ST-MOR04, BCDZ-MB03, PZ02, CO01, SS01, CEJO-DAM00, CPZ-CMA00, Sooryanarayana98, KRR-DAM96, Yushmanov87, SZ-CMJ04, CGH, Slater-JMPS88, SZ-CMJ03, MRRC-JDMSC06}, and arise in diverse areas including coin weighing problems \citep{Lindstrom64, ER63, SS-AMM63, KLT00, GN-AMM95}, network discovery and verification \citep{BEEHHMS-IEEE06}, robot navigation \citep{KRR-DAM96, SSH02}, connected joins in graphs \citep{ST-MOR04}, the Djokovi\'c-Winkler relation \citep{Caceres-etal}, and strategies for the Mastermind game \citep{Greenwell-JRM99, Goddard-JCMCC03, Goddard-JCMCC04, KLT00, Chvatal-Comb83}.

For positive integers $\beta$ and $D$, let \G\ be the class of connected graphs with metric dimension $\beta$ and diameter $D$. Consider the following two extremal questions:
\begin{itemize}
\item What is the minimum order of a graph in \G?
\item What is the maximum order of a graph in \G?
\end{itemize}

The first question was independently answered by \citet{Yushmanov87}, \citet{KRR-DAM96}, and \citet{CEJO-DAM00}, who proved that the minimum order of a graph in \G\ is $\beta+D$ (see \lemref{MinimumOrder}). Thus it is natural to consider the following problem:
\begin{itemize}
\item Characterise the graphs in \G\ with order $\beta+D$. 
\end{itemize}

Such a characterisation is simple for $\beta=1$. In particular, \citet{KRR-DAM96} and \citet{CEJO-DAM00} independently proved that paths $P_n$ (with $n\geq2$ vertices) are the only graphs with metric dimension $1$. Thus $\G[1,D]=\{P_{D+1}\}$. 

The characterisation is again simple at the other extreme with $D=1$. In particular, \citet{CEJO-DAM00} proved that the complete graph $K_n$ (with $n\geq1$ vertices) is the only graph with metric dimension $n-1$ (see \propref{DiameterOne}). Thus $\G[\beta,1]=\{K_{\beta+1}\}$. 

\citet{CEJO-DAM00} studied the case $D=2$, and obtained a non-trivial characterisation of graphs in \G[\beta,2] with order $\beta+2$ (see \propref{DiameterTwo}). 

The first contribution of this paper is to characterise the graphs in \G\ with order $\beta+D$ for all values of $\beta\geq1$ and $D\geq3$, thus completing the characterisation for all values of $D$. This result is stated and proved in \secref{MinOrder}.

We then study the second question above: What is the maximum order of a graph in \G? Previously, only a weak upper bound was known. In particular, \citet{KRR-DAM96} and \citet{CEJO-DAM00} independently proved that every graph in \G\ has at most $D^\beta+\beta$ vertices. This bound is tight only for $D\leq3$ or $\beta=1$. 

Our second contribution is to determine the (exact) maximum order of a graph in \G\ for all values of $D$ and $\beta$. This result is stated and proved in \secref{MaxOrder}.

\SECT{MinOrder}{Graphs with Minimum Order}

In this section we characterise the graphs in \G\ with minimum order. We start with an elementary lemma.

\LEM{Comp}{Let $S$ be a set of vertices in a connected graph $G$. Then $V(G)\setminus S$ resolves $G$ if and only if every pair of vertices in $S$ are resolved by some vertex not in $S$.}

\PROOF{If $v\in V(G)\setminus S$ and $w$ is any other vertex, then $v$ resolves $v$ and $w$. By assumption every pair of vertices in $S$ are resolved by some vertex in $V(G)\setminus S$.}



\lemref{Comp} enables the minimum order of a graph in \G\ to be easily determined.

\LEMMA{\citep{CEJO-DAM00,KRR-DAM96,Yushmanov87}}{MinimumOrder}{The minimum order of a graph in \G\ is $\beta+D$.}

\PROOF{First we prove that every graph $G\in\G$ has order at least $\beta+D$. Let $v_0,v_D$ be vertices such that $\dist(v_0,v_D)=D$. Let $P=(v_0,v_1,\dots,v_D)$ be a path of length $D$ in $G$. Then $v_0$ resolves $v_i,v_j$ for all distinct $i,j\in[1,D]$. Thus $V(G)\setminus\{v_1,\dots,v_D\}$ resolves $G$ by \lemref{Comp}. Hence $\beta\leq|V(G)|-D$ and $|V(G)|\geq\beta+D$.

It remains to construct a graph $G\in\G$ with order $\beta+D$. Let $G$ be the `broom' tree obtained by adding $\beta$ leaves adjacent to one endpoint of the path on $D$ vertices. Observe that $|V(G)|=\beta+D$ and $G$ has diameter $D$. It follows from Slater's formula \citep{Slater75} for the metric dimension of a
tree\footnote{Also see \citep{KRR-DAM96,CEJO-DAM00,HM-AC76} for proofs of Slater's formula.} that the $\beta$ leaves adjacent to one endpoint of the path are a metric basis of $G$. Hence $G\in\G$.
}

\SUBSECT{Twins}{Twin Vertices}

Let $u$ be a vertex of a graph $G$. The \emph{open neighborhood} of $u$ is $N(u):=\{v\in V(G):uv\in E(G)\}$, and the \emph{closed neighborhood} of $u$ is $N[u]:=N(u)\cup\{u\}$. Two distinct vertices $u,v$ are \emph{adjacent twins} if $N[u]=N[v]$, and \emph{non-adjacent twins} if $N(u)=N(v)$. Observe that if $u,v$ are adjacent twins then $uv\in E(G)$, and if $u,v$ are non-adjacent twins then $uv\not\in E(G)$; thus the names are justified\footnote{In the literature, adjacent twins are called \emph{true} twins, and non-adjacent twins are called \emph{false} twins. We prefer the more descriptive names, \emph{adjacent} and \emph{non-adjacent}.}. If $u,v$ are adjacent or non-adjacent twins, then $u,v$ are \emph{twins}. The next lemma follows from the definitions.

\LEM{TwinDistance}{If $u,v$ are twins in a connected graph $G$, then $\dist(u,x)=\dist(v,x)$ for every vertex $x\in V(G)\setminus\{u,v\}$.\qed}

\COR{totTwinMenysUn}{Suppose that $u,v$ are twins in a connected graph $G$ and $S$ resolves $G$. Then $u$ or $v$ is in $S$. Moreover, if $u\in S$ and $v\notin S$, then $(S\setminus\{u\})\cup\{v\}$ also resolves $G$.\qed}

\LEM{BasicTwin}{In a set $S$ of three vertices in a graph, it is not possible that two vertices in $S$ are adjacent twins, and two vertices in $S$ are non-adjacent twins.}

\PROOF{Suppose on the contrary that $u,v$ are adjacent twins and $v,w$ are non-adjacent twins. Since $u,v$ are twins and $v\not\sim w$, we have $u\not\sim w$. Similarly, since $v,w$ are twins and $u\sim v$, we have $u\sim w$. This is the desired contradiction.}

\LEM{twinsIguals}{Let $u,v,w$ be distinct vertices in a graph. If $u,v$ are twins and $v,w$ are twins, then $u,w$ are also twins.}

\PROOF{Suppose that $u,v$ are adjacent twins. That is, $N[u]=N[v]$. By \lemref{BasicTwin}, $v,w$ are adjacent twins. That is, $N[v]=N[w]$. Hence $N[u]=N[w]$. That is, $u,w$ are adjacent twins. By a similar argument, if 
$u,v$ are non-adjacent twins, then $v,w$ are non-adjacent twins and $u,w$ are non-adjacent twins.}



For a graph $G$, a set $T\subseteq V(G)$ is a \emph{twin-set} of $G$ if $v,w$ are twins in $G$ for every pair of distinct vertices $v,w\in T$. 

\LEM{TwinSetBasic}{If $T$ is a twin-set of a graph $G$, then either every pair of vertices in $T$ are adjacent twins, or every pair of vertices in $T$ are non-adjacent twins.}

\PROOF{Suppose on the contrary some pair of vertices $v,w\in T$ are adjacent twins, and some pair of vertices $x,y\in T$ are adjacent twins. If $v,x$ are adjacent twins then $\{v,x,y\}$ contradict \lemref{BasicTwin}. Otherwise 
$v,x$ are non-adjacent twins, in which case $\{v,w,x\}$ contradict \lemref{BasicTwin}.}

\LEM{DeleteTwin}{Let $T$ be a twin-set of a connected graph $G$ with $|T|\geq3$. Then $\beta(G)=\beta(G\setminus u)+1$ for every vertex $u\in T$.}

\PROOF{Let $u,v,w$ be distinct vertices in $T$. By \corref{totTwinMenysUn}, there is a metric basis $W$ of $G$ such that $u,v\in W$. Since $u$ has a twin in $G\setminus u$, for all $x,y\in V(G\setminus u)$ we have $\dist_{G}(x,y)=\dist_{G\setminus u}(x,y)$. In particular, $G\setminus u$ is connected. First we prove that $W\setminus\{u\} $ resolves $G\setminus u$. For all distinct vertices $x,y\in V(G\setminus u)$, there is a vertex $s\in W$ such that $\dist_G(x,s)\neq \dist_G(y,s)$. If $s\neq u$, then $s\in W\setminus\{u\}$ resolves the pair $x,y$. Otherwise, $v$ is a twin of $s=u$ and $\dist_{G\setminus u}(x,v)= \dist_{G}(x,v) =\dist_G(x,s)\neq\dist_G(y,s)=\dist_G(y,v)= 
\dist_{G\setminus u}(y,v)$. Consequently, $v\in W\setminus \{ u\}$ resolves the pair $x,y$. Now suppose that $W'$ is a resolving set of $G\setminus u$ such that $|W'|<|W|-1$. For all $x,y\in V(G\setminus u)$, there exists  a vertex  $s\in W'$ such that $\dist_{G\setminus u}(x,s)\neq \dist_{G\setminus u}(y,s)$. Then $W'\cup \{ u \}$ is a resolving set in $G$ of cardinality less than $|W|$, which contradicts the fact that $W$ is a resolving set of minimum cardinality.}

Note that it is necessary to assume that $|T|\geq3$ in \lemref{DeleteTwin}. For example, $\{x,z\}$ is a twin-set of the $3$-vertex path $P_3=(x,y,z)$, but $\beta(P_3)=\beta(P_3\setminus x)=1$.

\COR{quitarTwins}{Let $T$ be a twin-set of a connected graph $G$ with $|T|\geq3$. Then $\beta(G)=\beta(G\setminus S)+|S|$ for every subset $S\subset T$ with $|S|\leq|T|-2$.\qed}

\SUBSECT{TwinGraph}{The Twin Graph}

Let $G$ be a graph. Define a relation $\equiv$ on $V(G)$ by $u\equiv v$ if and only if $u=v$ or $u,v$ are twins. By \lemref{twinsIguals}, $\equiv$ is an equivalence relation. For each vertex $v\in V(G)$, let $v^*$ be the set of vertices of $G$ that are equivalent to $v$ under $\equiv$. Let $\{v_1^*,\dots,v_k^*\}$ be the partition of $V(G)$ induced by $\equiv$, where each $v_i$ is a representative of the set $v_i^*$. The \emph{twin graph} of $G$, denoted by $G^*$, is the graph with vertex set $V(G^*):=\{v_1^*,\dots,v_k^*\}$,  where $v_i^*v_j^*\in E(G^*)$ if and only if $v_iv_j\in E(G)$. The next lemma implies that this definition is independent of the choice of representatives.

\LEM{Reps}{Let $G^*$ be the twin graph of a graph $G$. Then two vertices $v^*$ and $w^*$ of $G^*$ are adjacent if and only if every vertex in $v^*$ is adjacent to every vertex in $w^*$ in $G$.}

\PROOF{Suppose on the contrary that some vertex in $v^*$ is adjacent to some vertex in $w^*$, and some vertex in $v^*$ is not adjacent to some vertex in $w^*$. Then $y\sim x\not\sim z$ for some vertices $x\in v^*$ and $y,z\in w^*$. Thus $y,z$ are not twins, which is the desired contradiction.}

Each vertex $v^*$ of $G^*$ is a maximal twin-set of $G$. By \lemref{TwinSetBasic}, $G[v^*]$ is a complete graph if the vertices of $v^*$ are adjacent twins, or $G[v^*]$ is a null graph if the vertices of $v^*$ are non-adjacent twins. So it makes sense to consider the following types of vertices in $G^*$. We say that $v^*\in V(G^*)$ is of \emph{type}:
\begin{itemize}
  \item $(1)$ if $|v^*| =1$,
  \item $(K)$ if $G[v^*]\isomorphic K_r$ and $r\geq2$,
  \item $(N)$ if $G[v^*]\isomorphic N_r$ and $r\geq2$, \\
  where $N_r$ is the \emph{null} graph with $r$ vertices and no edges.
\end{itemize}
A vertex of $G^*$ is of \emph{type} $(1K)$ if it is of type $(1)$ or $(K)$.
A vertex of $G^*$ is of \emph{type} $(1N)$ if it is of type $(1)$ or $(N)$.
A vertex of $G^*$ is of \emph{type} $(KN)$ if it is of type $(K)$ or $(N)$.

Observe that the graph $G$ is uniquely determined by $G^*$, and the type and cardinality of each vertex of $G^*$. In particular, if $v^*$ is adjacent to $w^*$ in $G^*$, then every vertex in $v^*$ is adjacent to every vertex in $w^*$ in $G$. 

We now show that the diameters of $G$ and $G^*$ are closely related.

\LEM{EqualDiameters}{Let $G\ne K_1$ be a connected graph. Then $\diam(G^*)\leq\diam(G)$. Moreover, $\diam(G^*)<\diam(G)$ if and only if $G^*\isomorphic K_n$ for some $n\geq1$. In particular, if $\diam(G)\geq3$ then $\diam(G)=\diam(G^*)$.}

\PROOF{If $v,w$ are adjacent twins in $G$, then $\dist_G(v,w)=1$ and $v^*=w^*$.
If $v,w$ are non-adjacent twins in $G$, then (since $G$ has no isolated vertices) $\dist_G(v,w)=2$ and $v^*=w^*$. If $v,w$ are not twins, then there is a shortest path between $v$ and $w$ that contains no pair of twins (otherwise there is a shorter path); thus \EQN{EqualDiameters}{\dist_G(v,w)=\dist_{G^*}(v^*,w^*).}

This implies that $\diam(G^*)\leq\diam(G)$. Moreover, if $\ecc_G(v)\geq3$ then $v$ is not a twin of every vertex $w$ for which $\dist_G(v,w)=\ecc_G(v)$; thus $\dist_G(v,w)=\dist_{G^*}(v^*,w^*)$ by \eqnref{EqualDiameters} and $\ecc_G(v)=\ecc_{G^*}(v^*)$. Hence if $\diam(G)\geq3$ then $\diam(G)=\diam(G^*)$. 

Now suppose that $\diam(G)>\diam(G^*)$.  Thus $\diam(G)\leq2$. If $\diam(G)=1$ then $G$ is a complete graph and $G^*\isomorphic K_1$, as claimed. Otherwise $\diam(G)=2$ and $\diam(G^*)\leq1$; thus $G^*\isomorphic K_n$ for some $n\geq1$, as claimed.

It remains to prove that $\diam(G^*)<\diam(G)$ whenever $G^*\isomorphic K_n$. In this case, $\diam(G^*)\leq 1$. So we are done if $\diam(G)\geq2$. Otherwise $\diam(G)\leq1$ and $G$ is also a complete graph. Thus $G^*\isomorphic K_1$ and $\diam(G^*)=0$. Since $G\neq K_1$, we have $\diam(G)=1>0=\diam(G^*)$, as desired.}

Note that graphs with $\diam(G^*)<\diam(G)$ include the complete multipartite graphs. 


\thmref{Characterisation} below characterizes the graphs in \G\ for $D\geq3$ in terms of the twin graph. \citet{CEJO-DAM00} characterized\footnote{To be more precise, \citet{CEJO-DAM00} characterised the graphs with $\beta(G)=n-2$. By \lemref{MinimumOrder}, if $\beta(G)=n-2$ then $G$ has diameter at most $2$. By \propref{DiameterOne}, if $G$ has diameter $1$ then $\beta(G)=n-1$. Thus if $\beta(G)=n-2$ then $G$ has diameter $2$.} the graphs in \G\ for $D\leq 2$. For consistency with \thmref{Characterisation}, we describe the characterisation by \citet{CEJO-DAM00} in terms of the twin graph.

\PROPOS{\citep{CEJO-DAM00}}{DiameterOne}{The following are equivalent for a  connected graph $G$ with $n$ vertices:
\begin{itemize}
\item $G$ has metric dimension $\beta(G)=n-1$,
\item $G\isomorphic K_n$,
\item $\diam(G)=1$,
\item the twin graph $G^*$ has one vertex, which is of type $(1K)$.
\end{itemize}}


\PROPOS{\citep{CEJO-DAM00}}{DiameterTwo}{The following are equivalent for a connected graph $G$ with $n\geq3$ vertices:
\begin{itemize}
\item $G$ has metric dimension $\beta(G)=n-2$,
\item $G$ has metric dimension $\beta(G)=n-2$ and diameter $\diam(G)=2$,
\item the twin graph $G^*$ of $G$ satisfies
\begin{itemize}
\item $G^*\isomorphic P_2$ with at least one vertex of type $(N)$, or
\item $G^*\isomorphic P_3$ with one leaf of type $(1)$, the other leaf of type $(1K)$,\\ and the degree-$2$ vertex of type $(1K)$.
\end{itemize}
\end{itemize}}

To describe our characterisation we introduce the following notation. Let $P_{D+1}=(u_0,u_1,\dots,u_D)$ be a path of length $D$. As illustrated in \figref{Figure1}(a), for $k\in[3,D-1]$ let $P_{D+1,k}$ be the graph obtained from $P_{D+1}$ by adding one vertex adjacent to $u_{k-1}$. As illustrated in \figref{Figure1}(b), for $k\in[2,D-1]$ let $P'_{D+1,k}$ be the graph obtained from $P_{D+1}$ by adding one vertex adjacent to $u_{k-1}$ and $u_k$. 

\Figure{Figure1}{\includegraphics{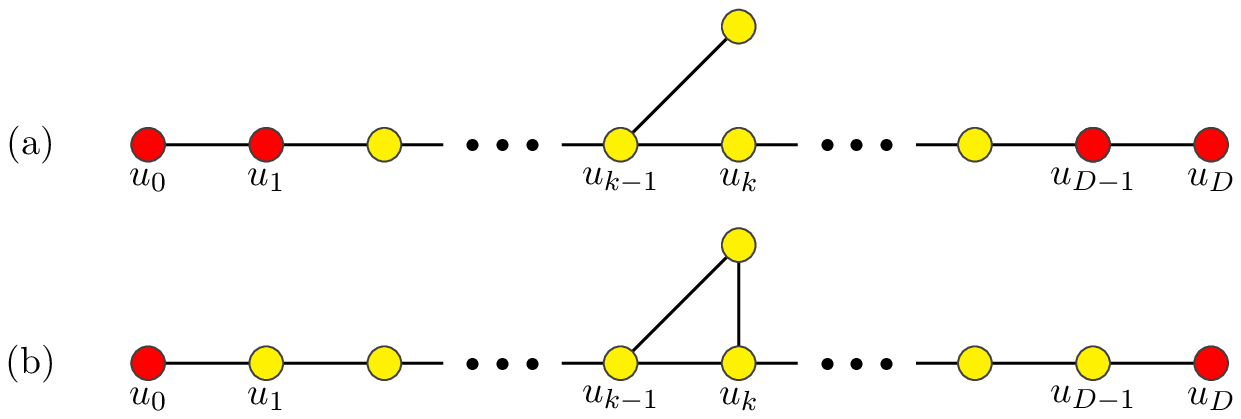}}{The graphs (a) $P_{D+1,k}$ and (b) $P'_{D+1,k}$.}

\THM{Characterisation}{Let $G$ be a connected graph of order $n$ and diameter $D\ge3$. Let $G^*$ be the twin graph of $G$. Let $\alpha(G^*)$ be the number
of vertices of $G^*$ of type $(K)$ or $(N)$. Then $\beta(G)=n-D$
if and only if $G^*$ is one of the following graphs:
\begin{enumerate}
  \item $G^*\isomorphic P_{D+1}$ and one of the following cases hold (see \figref{Figure2}):
\begin{enumerate}
  \item $\alpha(G^*)\le 1$;
  \item $\alpha(G^*)=2$, the two vertices of $G^*$ not of type $(1)$ are adjacent, and if one is a leaf of type $(K)$ then the other is also of type $(K)$;
  \item $\alpha(G^*)=2$, the two vertices of $G^*$ not of type $(1)$ are at distance $2$ and both are of type $(N)$; or
  \item $\alpha(G^*)=3$ and there is a vertex of type $(N)$ or $(K)$ adjacent to two vertices of type $(N)$.
\end{enumerate}
  \item $G^*\isomorphic P_{D+1,k}$ for some $k\in[3,D-1]$, 
the degree-$3$ vertex $u_{k-1}^*$ of $G^*$ is any type,
each neighbour of $u_{k-1}^*$ is type $(1N)$, 
and every other vertex is type $(1)$; see \figref{Figure3}.
  \item $G^*\isomorphic P'_{D+1,k}$  for some $k\in[2,D-1]$, 
the three vertices in the cycle are of type $(1K)$, 
and every other vertex is of type $(1)$; see \figref{Figure4}.
\end{enumerate}}

\Figure{Figure2}{\includegraphics{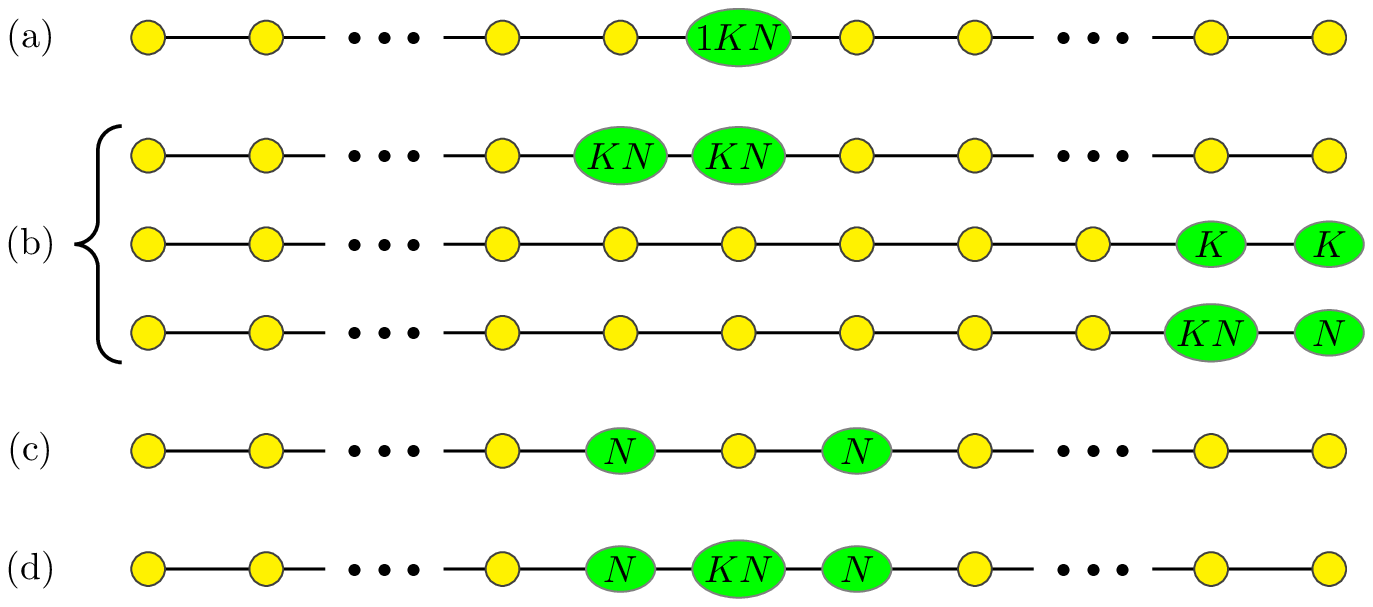}}{Cases (a)--(d) with $G^*\isomorphic P_{D+1}$ in \thmref{Characterisation}.}

\Figure{Figure3}{\includegraphics{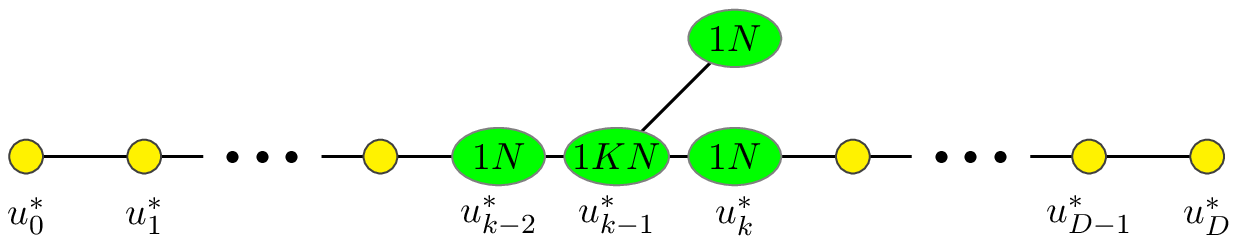}}{The case of $G^*\isomorphic P_{D+1,k}$ in \thmref{Characterisation}.}

\Figure{Figure4}{\includegraphics{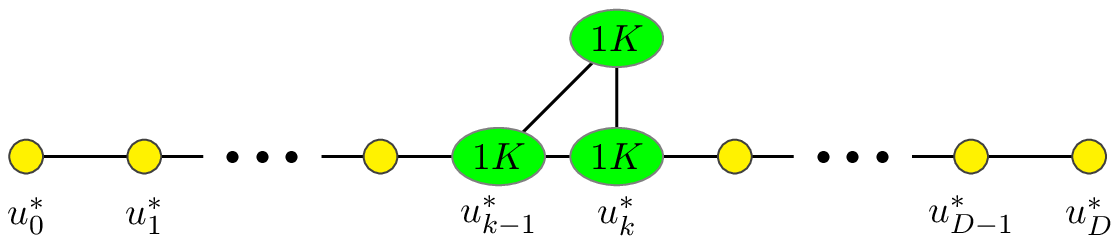}}{The case of $G^*\isomorphic P'_{D+1,k}$ in \thmref{Characterisation}.}

\SUBSECT{Necessity}{Proof of Necessity}

Throughout this section, $G$ is a graph of order $n$, diameter $D\ge 3$, and metric dimension $\beta(G)=n-D$. Let $G^*$ be the twin graph of $G$.

\LEM{excentricSenseTwins}{There exists a vertex $u_0$ in $G$ of eccentricity $D$ with no twin.}

\PROOF{Let $u_0$ and $u_D$ be vertices at distance $D$ in $G$. As illustrated in \figref{Figure5}, let $(u_0,u_1,\dots,u_D)$ be a shortest path between $u_0$ and $u_D$. Suppose on the contrary that both $u_0$ and $u_D$ have twins. Let $x$ be a twin of $u_0$ and $y$ be a twin of $u_D$. We claim that $\{x,y\}$ resolves $\{u_0,\dots,u_D\}$. Now $u_0\not\sim u_i$ 
for all $i\in[2,D]$, and thus $x\not\sim u_i$ (since $x,u_0$ are twins). Thus $\dist(x,u_i)=i$ for each $i\in[1,D]$. Hence $x$ resolves $u_i,u_j$ for all distinct $i,j\in[1,D]$. By symmetry, $\dist(y,u_i)=D-i$ for all $i\in[0,D-1]$, and $y$ resolves $u_i,u_j$ for all distinct $i,j\in[0,D-1]$. Thus $\{x,y\}$ resolves $\{u_0,\dots,u_D\}$, except for possibly the pair $u_0,u_D$. Now $\dist(x,u_0)\leq2$ and $\dist(x,u_D)=D$. Since $D\geq3$, $x$ resolves $u_0,u_D$. Thus $\{x,y\}$ resolves $\{u_0,\dots,u_D\}$. By \lemref{Comp}, $\beta(G)\leq n-(D+1)<n-D$, which is a contradiction. Thus $u_0$ or $u_D$ has no twin.}

\Figure{Figure5}{\includegraphics{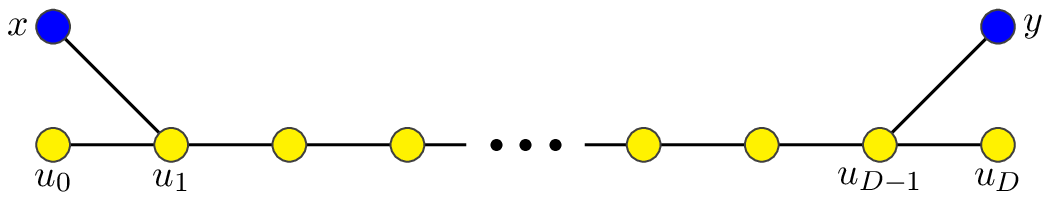}}{$\{x,y\}$ resolves $\{u_0,\dots,u_D\}$ in \lemref{excentricSenseTwins}.}

For the rest of the proof, fix a vertex $u_0$ of eccentricity $D$ in $G$ with no twin, which exists by \lemref{excentricSenseTwins}. Thus $u_0^*=\{u_0\}$ and $\ecc_{G^*}(u^*_0)=\ecc_G(u_0)=D$, which is also the diameter of $G^*$ by \lemref{EqualDiameters}. As illustrated in \figref{Figure6}, for each $i\in [0, D]$, let
\begin{align*}
A_i^*\,&:=\,\{v^*\in V(G^*)\,:\,\dist(u_0^*,v^*)=i\},\text{ and}\\
A_i\,&:=\,\{v\in V(G)\,:\,\dist(u_0,v)=i\}\,=\,\bigcup\{v^*:v^* \in A_i^*\}.
\end{align*}
Note that the last equality is true because $u_0$ has no twin and
$\dist(u_0,v)=\dist(u_0,w)$ if $v,w$ are twins. For all $i\in [0,D]$, we have $|A_i|\ge 1$ and $|A_i^*|\ge 1$. Moreover, $|A_0|=|A_0^*|=1$. Let $(u_0,u_1,\dots,u_D)$ be a path in $G$ such that $u_i\in A_i$ for each $i\in[0,D]$. Observe that if $v\in A_i$ is adjacent to $w\in A_j$ then $|i-j|\leq 1$. In particular, $(u_i,u_{i+1},\dots,u_j)$ is a shortest path between $u_i$ and $u_j$. 

\Figure{Figure6}{\includegraphics{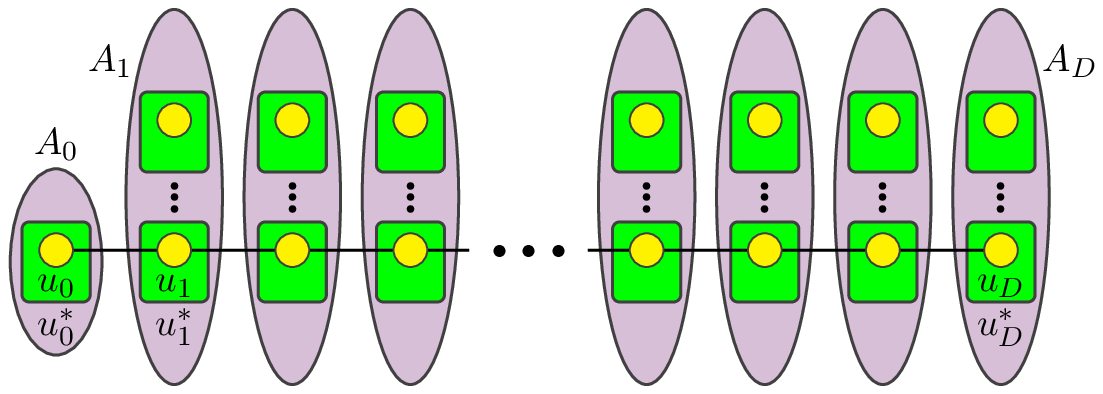}}{The sets $A_0,A_1,\dots,A_D$.}


\LEM{CompleteOrNull}{For each $k\in [1,D]$,
\begin{itemize}
  \item $G[A_k]$ is a complete graph or a null graph;
  \item $G^*[A_k^*]$ is a complete graph or a null graph, and all the vertices in $A_k^*$ are of type $(1K)$ in the first case, and of type $(1N)$ in the second case.
\end{itemize}}

\PROOF{Suppose that $G[A_k]$ is neither complete nor null for some $k\in[1,D]$. Thus there exist vertices $u,v,w\in A_k$ such that $u\sim v \not\sim w$, as illustrated in \figref{Figure7}\footnote{In \manyfigref{Figure7}{Figure23}, a solid line connects adjacent vertices, a dashed line connects non-adjacent vertices, and a coil connects vertices that may or may not be adjacent.}. Let $S:=(\{u_1,\dots,u_D\}\setminus\{u_k\})\cup\{u,w\}$. Every pair of vertices in $S$ are resolved by $u_0$, except for $u$ and $w$ which are resolved by $v$. Thus $\{u_0,v\}$ resolves $S$. By \lemref{Comp}, $\beta (G)\le
n-(D+1)<n-D$. This contradiction proves the first claim, which immediately implies the second claim.}

\Figure{Figure7}{\includegraphics{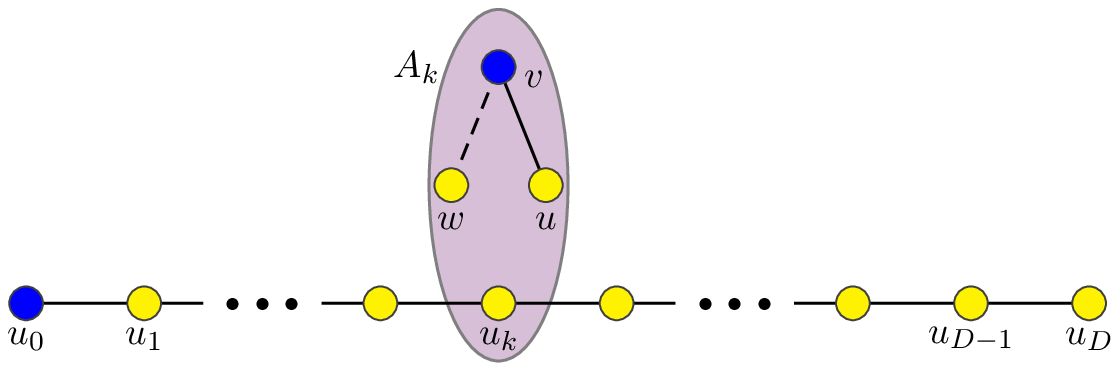}}{$\{u_0,v\}$ resolves  $(\{u_1,\dots,u_D\}\setminus\{u_k\})\cup\{u,w\}$ in \lemref{CompleteOrNull}.}


\LEM{BetweenLayers}{For each $k\in[1,D]$, if $|A_{k}|\geq2$ then
\begin{enumerate}
  \item[\textup{(a)}] $v\sim w$ for all vertices $v\in A_{k-1}$ and $w\in A_k$;
  \item[\textup{(b)}] $v^*\sim w^*$ for all vertices $v^*\in A_{k-1}^*$ and $w^*\in A_k^*$.
\end{enumerate}}

\PROOF{First we prove (a). Every vertex in $A_1$ is adjacent to $u_0$, which is the only vertex in $A_0$. Thus (a) is true for $k=1$. Now assume that $k\ge2$. Suppose on the contrary that $v\not\sim w$ for some $v\in A_{k-1}$ and $w\in A_k$. There exists a vertex $u\in A_{k-1}$ adjacent to $w$. As illustrated in \figref{Figure8}, if $w\neq u_k$ then $\{u_0,w\}$ resolves  $(\{u_1,\dots,u_D\}\setminus\{u_{k-1}\})\cup\{u,v\}$. 

\Figure{Figure8}{\includegraphics{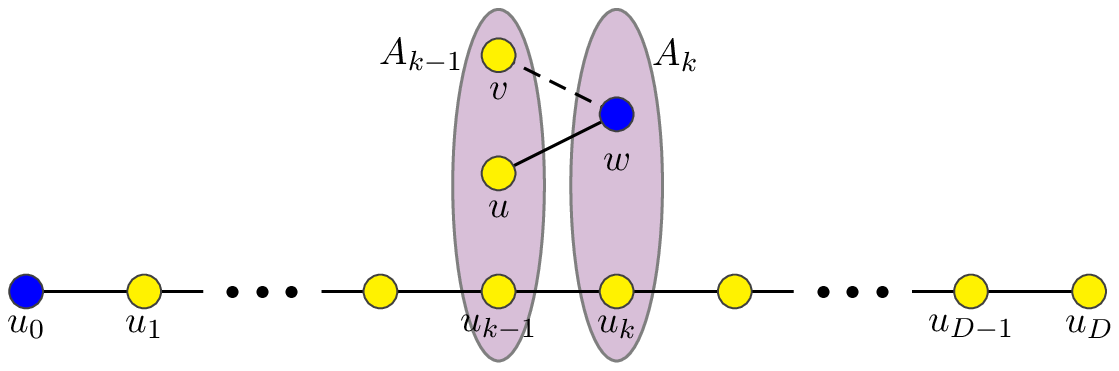}}{In \lemref{BetweenLayers}, $\{u_0,w\}$ resolves  $(\{u_1,\dots,u_D\}\setminus\{u_{k-1}\})\cup\{u,v\}$.}

As illustrated in \figref{Figure8b}, if $w=u_k$ then  $v\neq u_{k-1}$ and there exists a vertex $z\neq u_k$ in $A_k$, implying $\{u_0,u_k\}$ resolves $(\{u_1,\dots,u_D\}\setminus\{u_k\})\cup\{v,z\}$. In both cases, \lemref{Comp} implies that $\beta(G)\le n-D-1$. This contradiction proves (a), which immediately implies (b).}

\Figure{Figure8b}{\includegraphics{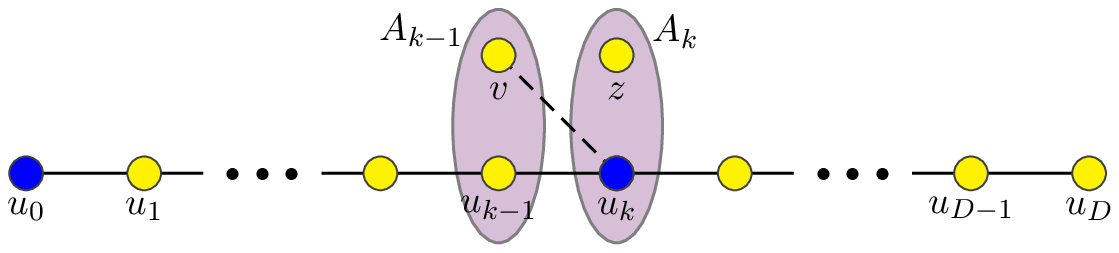}}{In \lemref{BetweenLayers}, $\{u_0,u_k\}$ resolves $(\{u_1,\dots,u_D\}\setminus\{u_k\})\cup\{v,z\}$.}


\LEM{BigAreClose}{If $|A_i|\geq2$ and $|A_j|\geq2$ then $|i-j|\leq 2$. Thus there are at most three distinct subsets $A_i,A_j,A_k$ each with cardinality at least $2$.}

\PROOF{As illustrated in \figref{Figure9}, suppose on the contrary that $|A_i|\ge 2$ and $|A_j|\ge 2$ for some $i,j\in[1,D]$ with $j\geq i+3$. Let $x\neq u_i$ be a vertex in $A_i$. Let $y\neq u_j$ be a vertex in $A_j$. We claim that $\{u_j,x\}$ resolves  $(\{u_0,\dots,u_D\}\setminus\{u_j\})\cup\{y\}$. 

By \lemref{BetweenLayers}, $u_{i-1}\sim x$ and $u_{j-1}\sim y$.
Observe that $\dist(u_j,y)\in \{ 1,2 \}$; $\dist(u_j,u_{j-h})=h$ for all $h\in[1,j]$; $\dist(u_j,u_{j+h})=h$ for all $h\in [1,D-j]$. Thus $u_j$ resolves $(\{u_0,\dots,u_D\}\setminus\{u_j\})\cup\{y\}$, except for the following pairs:
\begin{itemize}
\item  $u_{j-h},u_{j+h}$ whenever $1\leq h\leq j\leq D-h$;
\item  $y,u_{j-1}$ and $y,u_{j+1}$ if $\dist(y,u_j)=1$; and
\item   $y,u_{j-2}$ and $y,u_{j+2}$ if $\dist(y,u_j)=2$.
\end{itemize}

We claim that $x$ resolves each of these pairs. By \lemref{BetweenLayers}, there is a shortest path between $x$ and $u_{j-1}$ that passes through $u_{j-2}$. Let $r:=\dist(x,u_{j-2})$. Thus $\dist(x,u_{j-1})=r+1$, $\dist(x,y)=r+2$, $\dist(x,u_{j+1})=r+3$, and $\dist(x,u_{j+2})=r+4$. Thus $x$ resolves every pair of vertices in $\{u_{j-2},u_{j-1},y,u_{j+1},u_{j+2}\}$. It remains to prove that $x$ resolves $u_{j-h},u_{j+h}$ whenever $3\leq h\leq j\leq D-h$. Observe that $\dist(x,u_{j+h})\geq j+h-i$. If $j-h\geq i$ then, since $(x,u_{i-1},\dots,u_{j-h})$ is a path,
\eqn{\dist(x,u_{j-h})\leq j-h-i+2<j+h-i\leq\dist(x,u_{j+h}).} 
Otherwise $j-h\le i-1$, implying 
\eqn{\dist(x,u_{j-h})=i-(j-h)<j+h-i\leq\dist(x,u_{j+h}).} 
In each case $\dist(x,u_{j-h})<\dist(x,u_{j+h})$. Thus $x$ resolves $u_{j-h},u_{j+h}$.

Hence $\{u_j,x\}$ resolves  $(\{u_0,\dots,u_D\}\setminus\{u_j\})\cup\{y\}$. By \lemref{Comp}, $\beta(G)\le n-D-1$ which is the desired contradiction.}


\Figure{Figure9}{\includegraphics{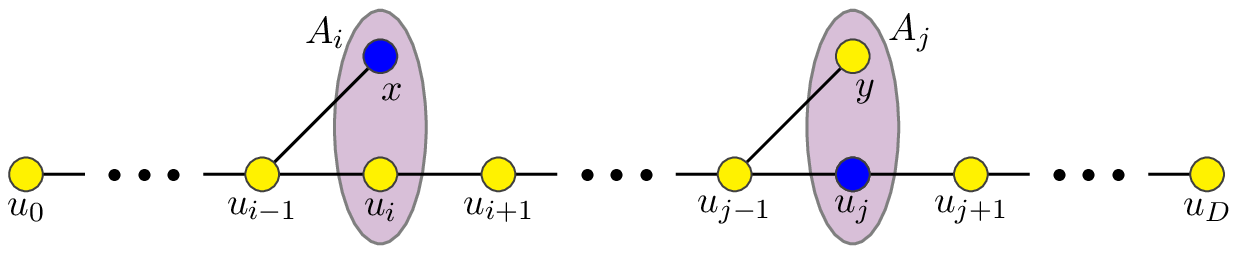}}{$\{u_j,x\}$ resolves  $(\{u_0,\dots,u_D\}\setminus\{u_j\})\cup\{y\}$ in \lemref{BigAreClose}.}

\LEM{A1AD}{$|A_1^*|=1$ and $|A_D^*|=1$.}

\PROOF{Consider a vertex $v\in A_1$. Then $v\sim u_0$ and every other neighbour of $v$ is in $A_1\cup A_2$. By \lemref{CompleteOrNull}, $G[A_1]$ is complete or null. If every vertex in $A_1$ is adjacent to every vertex in $A_2$, then $A_1$ is a twin-set, and $|A_1^*|=1$ as desired. 

Now assume that some vertex $v\in A_1$ is not adjacent to some vertex in $A_2$. By \lemref{BetweenLayers}, the only vertex in $A_2$ is $u_2$, and $v\not\sim u_2$. If $G[A_1]$ is null then $\ecc(v)>D$, and if $G[A_1]$ is complete then $v$ and $u_0$ are twins. In both cases we have a contradiction. 

If $|A_D|=1$ then $|A_D^*|=1$. Now assume that $|A_D|\geq2$. The neighbourhood of every vertex in $A_D$ is contained in $A_{D-1}\cup A_D$. By \lemref{BetweenLayers}, every vertex in $A_D$ is adjacent to every vertex in $A_{D-1}$. By \lemref{CompleteOrNull}, $G[A_D]$ is complete or null. Thus $A_D$ is a twin-set, implying $|A_D^*|=1$.}


\LEM{NewLemma}{For each $k\in[1,D-1]$, distinct vertices $v,w\in A_k$ are twins if and only if they have the same neighbourhood in $A_{k+1}$. }

\PROOF{The neighbourhood of both $v$ and $w$ is contained in $A_{k-1}\cup A_k\cup A_{k+1}$. By \lemref{BetweenLayers}, both $v$ and $w$ are adjacent to every vertex in $A_{k-1}$. By \lemref{CompleteOrNull}, $G[A_k]$ is complete or null. Thus $v$ and $w$ are twins if and only if they have the same neighbourhood in $A_{k+1}$.}


\LEM{cardinalAk}{For each $k\in[2,D]$,
\begin{enumerate}
  \item[\textup{(a)}] if $|A_k|\geq2$ then $|A_{k-1}^*|=1$;
  \item[\textup{(b)}] if $|A_k|=1$ then $|A_{k-1}^*|\le 2$.
\end{enumerate}}

\PROOF{Suppose that $|A_k|\geq2$. If $|A_{k-1}|=1$ then $|A_{k-1}^*|=1$ as desired. Now assume that $|A_{k-1}|\geq2$. Thus $A_{k-1}$ is a twin-set by \lemref{NewLemma}, implying $|A_{k-1}^*|=1$. Now suppose that $|A_{k}|=1$. If $|A_{k-1}|=1$ then $|A_{k-1}^*|=1$ and we are done. So assume that $|A_{k-1}|\geq2$. By \lemref{NewLemma}, the set of vertices in $A_{k-1}$ that are adjacent to the unique vertex in $A_{k}$ is a maximal twin-set, and the set of vertices in $A_{k-1}$ that are not adjacent to the unique vertex in $A_{k}$ is a maximal twin-set (if it is not empty). Therefore $|A_{k-1}^*|\leq2$.}


\LEM{LayerAtMostTwoBags}{For each $k\in[1,D]$, we have $|A_k^*|\leq 2$. Moreover, there are at most three values of $k$ for which $|A_k^*|=2$. Furthermore, if $|A_i^*|=2$ and $|A_j^*|=2$ then $|i-j|\leq2$.}

\PROOF{\lemref{A1AD} proves the result for $k=D$. Now assume that $k\in[1,D-1]$. Suppose on the contrary that $|A_k^*|\geq3$ for some $k\in[1,D]$. 
By the contrapositive of \lemref{cardinalAk}(a), $|A_{k+1}|=1$. By \lemref{cardinalAk}(b), $|A_k^*|\leq2$, which is the desired contradiction.
The remaining claims follow immediately from \lemref{BigAreClose}.}


\LEM{SmallAfterBig}{Suppose that $|A_k^*|=2$ for some $k\in[2,D-1]$. Then $|A_{k+1}|=|A_{k+1}^*|=1$, and exactly one of the two vertices of $A_k^*$ is adjacent to the only vertex of $A_{k+1}^*$. Moreover, if $k\le D-2$ then $|A_{k+2}|=|A^*_{k+2}|=1$.}

\PROOF{By the contrapositive of \lemref{cardinalAk}(a), $|A_{k+1}|=|A_{k+1}^*|=1$. By \lemref{NewLemma}, exactly one vertex in $A_k^*$ is adjacent to the vertex in $A_{k+1}^*$. Now suppose that $k\leq D-2$ but $|A_{k+2}|\geq2$. As illustrated in \figref{Figure12}, let $x\neq u_{k+2}$ be a vertex in $A_{k+2}$. Let $y\neq u_k$ be a vertex in $A_{k}$, such that $y,u_k$ are not twins, that is, $y\not\sim u_{k+1}$. By \lemref{BetweenLayers}, $u_{k-1}\sim y$ and $u_{k+1}\sim x$. Thus $\{x,u_0\}$ resolves $\{u_1,\dots ,u_D,y\}$. By \lemref{Comp}, $\beta(G)\leq n-D-1$, which is a contradiction.
Hence $|A_{k+2}|=1$, implying $|A^*_{k+2}|=1$. }

\Figure{Figure12}{\includegraphics{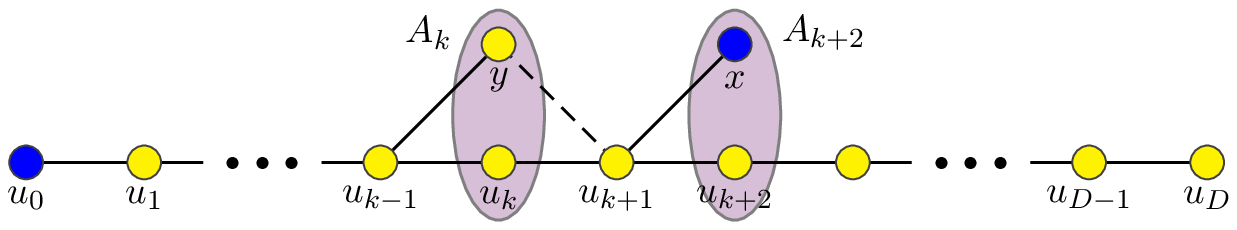}}{$\{x,u_0\}$ resolves $\{u_1,\dots,u_D,y\}$ in \lemref{SmallAfterBig}.}


We now prove that the structure of the graph $G^*$ is as claimed in \thmref{Characterisation}.

\LEM{ThreeTypes}{Either $G^*\isomorphic P_{D+1}$, $G^*\isomorphic P_{D+1,k}$ for some $k\in[3,D-1]$, or $G^*\isomorphic P'_{D+1,k}$ for some $k\in[2,D-1]$.}

\PROOF{By \lemref{LayerAtMostTwoBags} each set $A_k^*$ contains at most two vertices of $G^*$. \threelemref{A1AD}{BigAreClose}{SmallAfterBig} imply that $|A_k^*|=2$ for at most one $k\in[0,D]$. If $|A_k^*|=1$ for every $k\in[0,D]$ then $G^*\isomorphic P_{D+1}$ as desired. 

Now assume that $|A_k^*|=2$ for exactly one $k\in[0,D]$. By \lemref{A1AD}, $k\in[2,D-1]$. Let $w^*$ be the vertex in $A_k^*$ besides $u_k^*$. Then $w^*\sim u_{k-1}^*$ by \lemref{BetweenLayers}. If $w^*\sim u_k^*$ then $G^*\isomorphic P'_{D+1,k}$. Otherwise $w^*\not\sim u_k^*$. Then $G^*\isomorphic P_{D+1,k}$. It remains to prove that in this case $k\ne 2$. 

Suppose on the contrary that $G^*\isomorphic P_{D+1,k}$ and $k=2$. Thus $|A_2^*|=2$. Say $A_2^*=\{u_2^*,w^*\}$, where $u_2^*\not\sim w^*$. By \lemref{SmallAfterBig}, $|A_3^*|=1$. Thus $A_3^*=\{u_3^*\}$. Since $u_2^*\sim u_3^*$, by \lemref{NewLemma}, $w^*\not\sim u_3^*$. Thus $u_1^*$ is the only neighbour of $w^*$. Hence every vertex in $w^*$ is a twin of $u_0$, which contradicts the fact that $u_0$ has no twin. Thus $k\ne 2$ if $G^*\isomorphic P_{D+1,k}$.}




We now prove restrictions about the type of the vertices in $G^*$. To start with, \lemref{BigAreClose} implies:

\COR{cami}{If $G^*\isomorphic P_{D+1}$ then $\alpha(G^*)\le 3$
and the distance between every pair of vertices not of type $(1)$ is
at most $2$.}


\LEM{cami2junts}{Suppose that $G^*\isomorphic P_{D+1}$ and $\alpha(G^*)= 2$. If the two vertices of $G^*$ not of type $(1)$ are adjacent, and one of them is a leaf of type $(K)$, then the other is also of type $(K)$.}

\PROOF{As illustrated in \figref{Figure14}, let $x$ and $y$ be twins of $u_{D-1}$ and $u_D$ respectively. By assumption $G[A_D]$ is a complete graph. Suppose on the contrary that $G[A_{D-1}]$ is a null graph. By \lemref{BetweenLayers}, every vertex in $A_D$ is adjacent to every vertex in $A_{D-1}$. Thus $y$ resolves $\{u_0,\dots,u_D\}$, except for the pair $u_{D-1},u_D$, which is resolved by $x$. Thus $\{x,y\}$ resolves $\{u_0,\dots,u_D\}$. By \lemref{Comp}, $\beta(G)\leq n-D-1$, which is a contradiction. Thus $G[A_{D-1}]$ is a complete graph. }

\Figure{Figure14}{\includegraphics{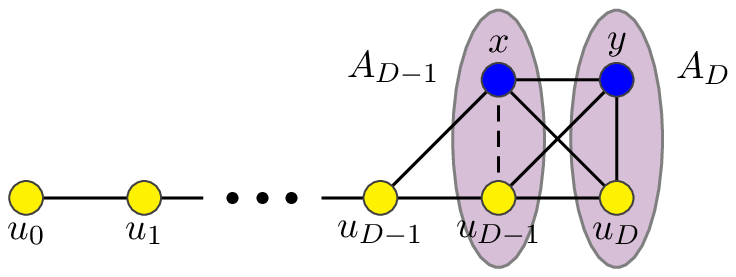}}{$\{x,y\}$ resolves $\{u_0,\dots,u_D\}$ in \lemref{cami2junts}.}

\LEM{cami2separats}{Suppose that $G^*\isomorphic P_{D+1}$ and for some $k\in[2,D-1]$, the vertices $u_{k-1}^*$ and $u_{k+1}^*$ of $G^*$ are both not of type (1).  Then $u_{k-1}^*$ and $u_{k+1}^*$ are both of type $(N)$.}

\PROOF{Let $x$ and $y$ be twins of $u_{k-1}$ and $u_{k+1}$ respectively. Suppose on the contrary that one of $u_{k-1}^*$ and  $u_{k+1}^*$ is of type $(K)$. Without loss of generality $u_{k-1}^*$ is of type $(K)$, as illustrated in \figref{Figure15}. Thus $u_{k-1}\sim x$. We claim that $\{x,y\}$ resolves $\{u_0,u_1,\dots,u_D\}$. 

Observe that $x$ resolves every pair of vertices of $\{u_0,u_1,\dots,u_D\}$ except for:

\begin{itemize}
\item each pair of vertices in $\{ u_{k-2},u_{k-1},u_k\}$, which are all resolved by $y$ since $d(y,u_k)=1$, $d(y,u_{k-1})=2$, and $d(y,u_{k-2})=3$; and

\item the pairs $\{u_{k-j},u_{k+j-2}:j\in[3,\min\{k,D+2-k\}]\}$, which are all resolved by $y$ since $d(y,u_{k-j})=j+1$, and 
\begin{equation*}
d(y,u_{k+j-2})=
\begin{cases}
j-2&\text{ if }j\geq4,\\
1\text{ or }2&\text{ if }j=3.
\end{cases}
\end{equation*}
\end{itemize}

Hence $\{x,y\}$ resolves $\{u_0,u_1,\dots,u_D\}$. Thus \lemref{Comp} implies $\beta(G)\le n-D-1$, which is the desired contradiction. Hence $u_{k-1}^*$ and  $u_{k+1}^*$ are both of type $(N)$.}


\Figure{Figure15}{\includegraphics{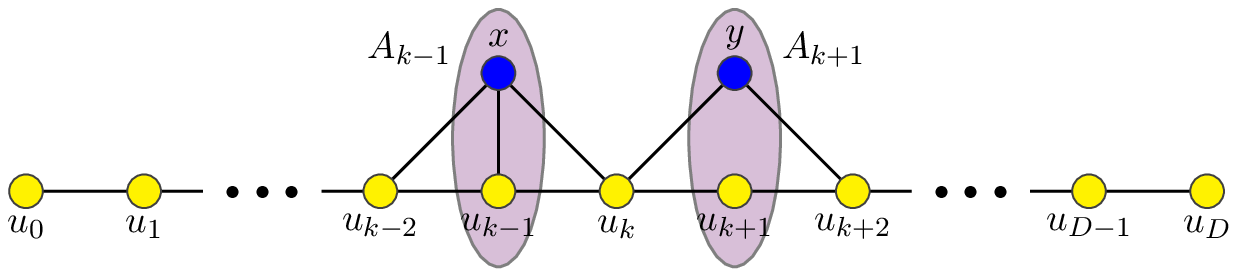}}{$\{x,y\}$ resolves
$\{u_0,\dots,u_D\}$ in \lemref{cami2separats}.}

\corref{cami} and \twolemref{cami2junts}{cami2separats} prove the necessity of the conditions in \thmref{Characterisation} when $G^*\isomorphic P_{D+1}$.


%



\LEM{camiMesFulla}{Suppose that $G^*\isomorphic P_{D+1,k}$ for some $k\in[3,D-1]$, where $A_k^*=\{u_k^*,w^*\}$ and $w^*\sim u^*_{k-1}$. Then $u_{k-2}^*$, $u^*_{k}$ and $w^*$ are type $(1N)$, $u_{k-1}^*$ is any type, and every other vertex is type $(1)$.}


\PROOF{Since $u^*_{k}\not\sim w^*$, \lemref{CompleteOrNull} implies that $u^*_{k}$ and $w^*$ are both type $(1N)$. By \twolemref{BigAreClose}{SmallAfterBig}, the remaining vertices are of type $(1)$ except, possibly $u^*_{k-2}$ and $u^*_{k-1}$. Suppose that $u^*_{k-2}$ is of type $(K)$, as illustrated in \figref{Figure16}. Let $x$ be a twin of $u_{k-2}$. Then $x\sim u_{k-1}$. We claim that $\{x,w\}$ resolves $\{u_0,u_1,\dots,u_D\}$. 

Observe that $x$ resolves every pair of vertices in $\{u_0,u_1,\dots,u_D\}$ except for:

\begin{itemize}
\item each pair of vertices in $\{ u_{k-3},u_{k-2},u_{k-1}\}$, which are all resolved by $w$ since $d(w,u_{k-1})=1$, $d(w,u_{k-2})=2$, and $d(w,u_{k-3})=3$; and

\item the pairs $\{u_{k-2-j},u_{k-2+j}:j\in[2,\min\{k-2, D-k+2\}]\}$, which are all resolved by $w$ since $d(w,u_{k-2-j})=j+2$ and $d(w,u_{k-2+j})=j$.
\end{itemize}

Thus $\{x,w\}$ resolves $\{u_0,u_1,\dots,u_D\}$. Hence \lemref{Comp} implies that $\beta(G)\leq n-D-1$, which is the desired contradiction.}

\Figure{Figure16}{\includegraphics{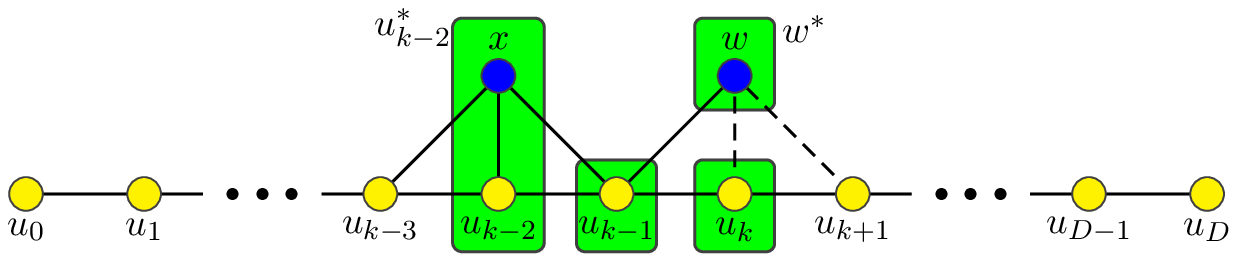}}{$\{x,w\}$ resolves  $\{u_0,\dots,u_D\}$ in \lemref{camiMesFulla}.}


\LEM{camiMesTriangle}{Suppose that $G^*\isomorphic P'_{D+1,k}$ for some $k\in[2,D-1]$, where $A_k^*=\{u_k^*,w^*\}$ and $u^*_{k-1}\sim w^*\sim u^*_k$. Then $u^*_{k-1}$, $u^*_{k}$ and $w^*$ are type $(1K)$, and every other vertex is type $(1)$.}

\PROOF{Since $u^*_{k}\sim w^*$, \lemref{CompleteOrNull} implies that $u^*_{k}$ and $w^*$ are type $(1K)$. By \twolemref{BigAreClose}{SmallAfterBig}, the remaining vertices are type $(1)$ except possibly $u^*_{k-2}$ and $u^*_{k-1}$. 

Suppose on the contrary that $u^*_{k-2}$ is type $(K)$ or $(N)$, as illustrated in \figref{Figure17}. Let $x$ be a twin of $u_{k-2}$. 
We claim that $\{x,w\}$ resolves $\{u_0,u_1,\dots,u_D\}$. 
Observe that $w$ resolves every pair of vertices in $\{u_0,u_1,\dots ,u_D\}$,
except for pairs 
\begin{equation*}
\{u_{k-1-j},u_{k+j}:j\in[0,\min\{k-1,D-k\}]\}.
\end{equation*}
These pairs are all resolved by $x$ since $d(x,u_{k+j})=j+2$ and 
\begin{equation*}
d(x,u_{k-1-j})=
  \begin{cases}
   j-1 & \text{ if } j\ge 2, \\
   1 \text{ or }2 & \text{ if } j=1, \\
   1 & \text{ if }j=0.
  \end{cases}
  \end{equation*}
Thus $\{x,w\}$ resolves $\{u_0,u_1,\dots,u_D\}$.

\Figure{Figure17}{\includegraphics{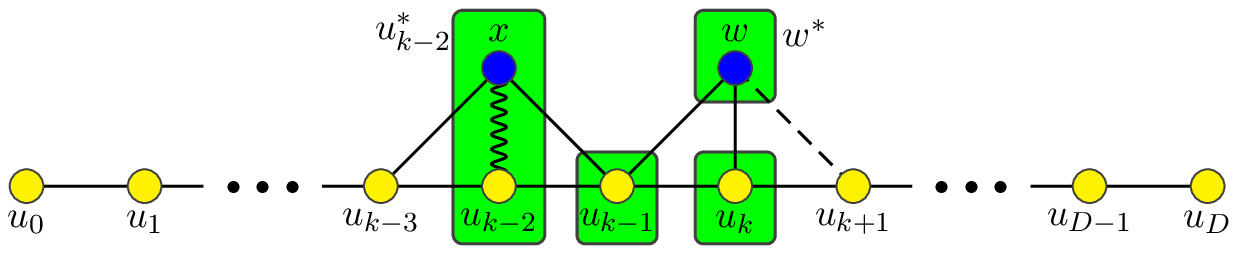}}{$\{x,w\}$ resolves $\{u_0,u_1,\dots,u_D\}$ in \lemref{camiMesTriangle}.}

Suppose on the contrary that $u^*_{k-1}$ is type $(N)$, as illustrated
in \figref{Figure17b}. Let $y$ be a twin of $u_{k-1}$. 
We claim that $\{y,w\}$ resolves $\{u_0,u_1,\dots,u_D\}$. 
Observe that $w$ resolves every pair of vertices in $\{u_0,u_1,\dots,u_D\}$,
except for pairs 
\begin{equation*}
\{u_{k-1-j},u_{k+j}:j\in[0,\min\{k-1,D-k\}]\}.
\end{equation*}
These pairs are all resolved by $y$ since $d(y,u_{k+j})=j+1$ and 
\begin{equation*}
d(y,u_{k-1-j})=
  \begin{cases}
  j	& \text{ if } j\geq 1, \\
  2	& \text{ if } j=0.
  \end{cases}
\end{equation*}
Thus $\{y,w\}$ resolves $\{u_0,u_1,\dots,u_D\}$. 

\Figure{Figure17b}{\includegraphics{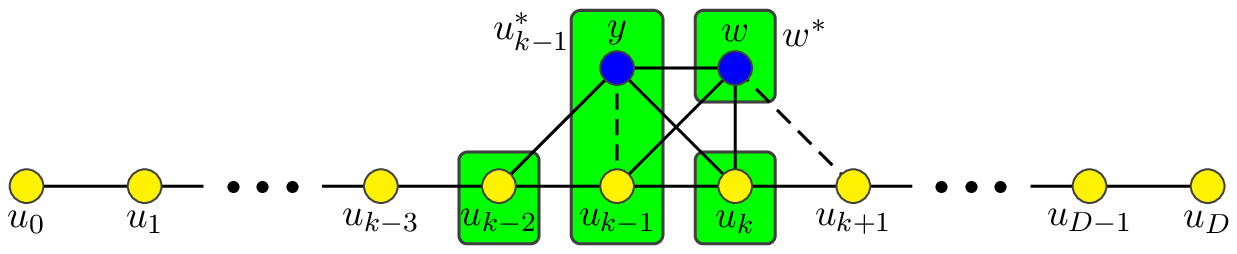}}{$\{y,w\}$ resolves $\{u_0,u_1,\dots,u_D\}$ in \lemref{camiMesTriangle}.}

By \lemref{Comp}, in each case $\beta(G)\leq n-D-1$, which is the desired contradiction.}

Observe that \twolemref{camiMesFulla}{camiMesTriangle} imply the necessity of the conditions in \thmref{Characterisation} when $G^*\isomorphic P_{D+1,k}$ or $G^*\isomorphic P'_{D+1,k}$. This completes the proof of the necessity of the conditions in \thmref{Characterisation}. 


\SUBSECT{Sufficiency}{Proof of Sufficiency}

Let $G$ be a graph with $n$ vertices and $\diam(G)\geq3$. Let $T$ be a twin-set of cardinality $r\geq 3$ in $G$. Let $G'$ be the graph obtained from $G$ by deleting all but two of the vertices in $T$. As in \lemref{EqualDiameters}, $\diam(G')=\diam(G)$. Say $G'$ has order $n'$. Then by \corref{quitarTwins}, $\beta(G')=\beta(G)-(r-2)$. Since $n'=n-(r-2)$, we have that $\beta(G)=n-D$ if and only if $\beta(G')=n'-D$. Thus it suffices to prove the sufficiency in \thmref{Characterisation} for graphs $G$ whose maximal twin-sets have at most two vertices. We assume in the remainder of this section that every twin-set in $G$ has at most two vertices.

Suppose that the twin graph $G^*$ of $G$ is one of the graphs stated in \thmref{Characterisation}. We need to prove that $\beta(G)=n-D$.  Since $\beta(G)\leq n-D$ by \lemref{MinimumOrder}, it suffices to prove that every subset of $n-D-1$ vertices of $G$ is not a resolving set. By \corref{totTwinMenysUn}, every resolving set contains at least one vertex in each twin-set of cardinality $2$. Observe also that, since $\alpha(G^*)$ is the number of vertices of $G^*$ not of type $(1)$, we have that $\alpha(G^*)=n-|V(G^*)|$.

\textbf{Case 1.} $G^*\isomorphic P_{D+1}$ with vertices $u_0^*\sim u^*_1\sim\dots\sim u^*_D$: We now prove that for each subcase stated in \thmref{Characterisation} every set of $n-D-1=n-|V(G^*)|=\alpha(G^*)$ vertices of $G$ does not resolve $G$. Suppose on the contrary that $W$ is a resolving set of $G$ of cardinality $\alpha(G^*)$. 

\textbf{Case 1(a).} $\alpha(G^*)\leq1$: We need at least one vertex to resolve a graph $G$ of order $n\ge 2$. So $\alpha(G^*)=1$. Thus $G$ is not a path, but \citet{KRR-DAM96} and \citet{CEJO-DAM00} independently proved that every graph with metric dimension $1$ is a path, which is a contradiction.

\textbf{Case 1(b)(i).} $\alpha(G^*)=2$, and  $u_k^*,u_{k+1}^*$ are not of type $(1)$  for some $k\in[1,D-2]$:\\ As illustrated in \figref{Figure18}, consider vertices $x\neq u_k$ in $u^*_k$, and $y\neq u_{k+1}$ in $u^*_{k+1}$. By \corref{totTwinMenysUn}, we may assume that $W=\{x,y\}$. 

Suppose that $u^*_k$ is type $(N)$. Then $x\not\sim u_k$, implying $\dist(x,u_k)=\dist(x,u_{k+2})=2$ and  $\dist(y,u_k)=\dist(y,u_{k+2})=1$. Thus neither $x$ nor $y$ resolve $u_k,u_{k+2}$.

Suppose that $u^*_{k+1}$ is type $(N)$. Then $y\not\sim u_{k+1}$, implying 
$\dist(x,u_{k-1})=\dist(x,u_{k+1})=1$ and  $\dist(y,u_{k-1})=\dist(y,u_{k+1})=2$. Thus neither $x$ nor $y$ resolve $u_{k-1},u_{k+1}$.

Suppose that $u^*_k$ and $u^*_{k+1}$ are both type $(K)$. Then $x\sim u_k$ and $y\sim u_{k+1}$, implying $\dist(x,u_{k})=\dist(x,u_{k+1})=1$ and  $\dist(y,u_{k})=\dist(y,u_{k+1})=1$. Thus neither $x$ nor $y$ resolve $u_k,u_{k+1}$.

In each case we have a contradiction.

\Figure{Figure18}{\includegraphics{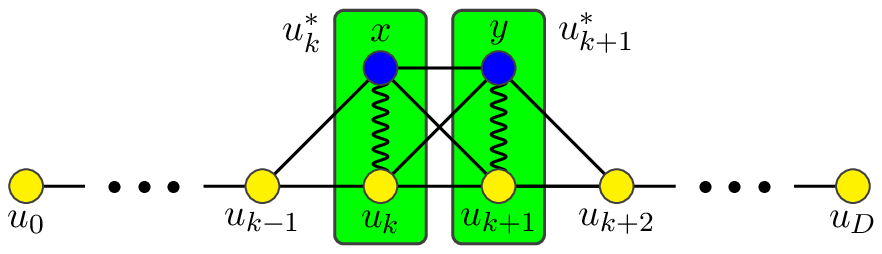}}{In Case 1(b)(i).}

\textbf{Case 1(b)(ii).} $\alpha (G^*)=2$, $u_{D-1}^*$ is not type (1), and $u_{D}^*$ is not type $(1)$: As illustrated in \figref{Figure19}, consider $x\neq u_{D-1}$ in $u^*_{D-1}$ and  $y\neq u_{D}$ in $u^*_{D}$. By \corref{totTwinMenysUn}, we may assume that $W=\{x,y\}$. 

First suppose that $u_D^*$ is of type $(N)$. Then $y\not\sim u_D$, 
implying $\dist(x,u_{D-2})=\dist(x,u_{D})=1$ and  $\dist(y,u_{D-2})=\dist(y,u_{D})=2$. Thus neither $x$ nor $y$ resolve $u_{D-2},u_{D}$, which is a contradiction.

Suppose that $u^*_D$ and $u^*_{D-1}$ are both type $(K)$. Then $x\sim u_{D-1}$ and $y\sim u_D$, implying $\dist(x,u_{D-1})=\dist(x,u_{D})=1$ and  $\dist(y,u_{D-1})=\dist(y,u_{D})=1$. Thus neither $x$ nor $y$ resolve $u_{D-1},u_D$, which is a contradiction.

\Figure{Figure19}{\includegraphics{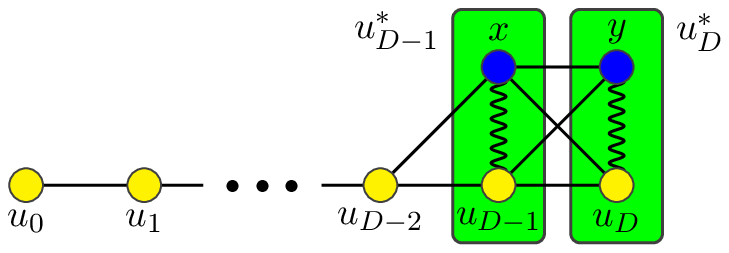}}{In Case 1(b)(ii).}


\textbf{Case 1(c).} 
$\alpha (G^*)=2$ and  $u_{k-1}^*$ is type ($N)$, and $u_{k+1}^*$ is type $(N)$ for some $k\in[2,D-1]$: As illustrated in \figref{Figure20}, consider $x\neq u_{k-1}$ in $u^*_{k-1}$ and $y\neq u_{k+1}$ in $u^*_{k+1}$. By \corref{totTwinMenysUn}, we may assume that $W=\{x,y\}$. Since  $x\not\sim u_{k-1}$ and $y\not\sim u_{k+1}$, we have $\dist(x,u_{k-1})=\dist(x,u_{k+1})=2$ and  $\dist(y,u_{k-1})=\dist(y,u_{k+1})=2$. Thus neither $x$ nor $y$ resolve $u_{k-1},u_{k+1}$, which is a contradiction.

\Figure{Figure20}{\includegraphics{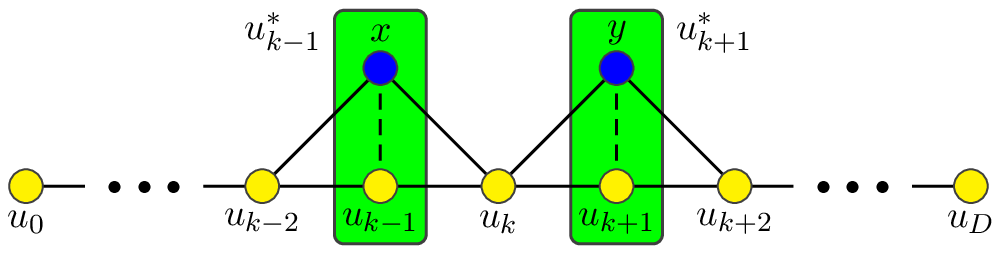}}{$\{x,y\}$ does not resolve $u_{k-1},u_{k+1}$ in Case 1(c).}


\textbf{Case 1(d).} $\alpha(G^*)=3$, $u_{k-1}^*$ is type $(N)$, $u_k^*$ is type $(K)$ or $(N)$, and $u_{k+1}^*$ is type $(N)$ for some $k\in [2,D-1]$: 
As illustrated in \figref{Figure21}, consider $x\neq u_{k-1}$ in $u^*_{k-1}$,  $y\neq u_{k}$ in $u^*_{k}$, and $z\neq u_{k+1}$ in $u^*_{k+1}$. By \corref{totTwinMenysUn}, we may assume that that $W=\{x,y,z\}$. Now $x\not\sim u_{k-1}$ and $z\not\sim u_{k+1}$. Thus $\dist(x,u_{k-1})=\dist(x,u_{k+1})=2$, $\dist(y,u_{k-1})=\dist(y,u_{k+1})=1$, and $\dist(z,u_{k-1})=\dist(z,u_{k+1})=2$. Thus $\{x,y,z\}$ does not resolve $u_{k-1},u_{k+1}$, which is a contradiction.

\Figure{Figure21}{\includegraphics{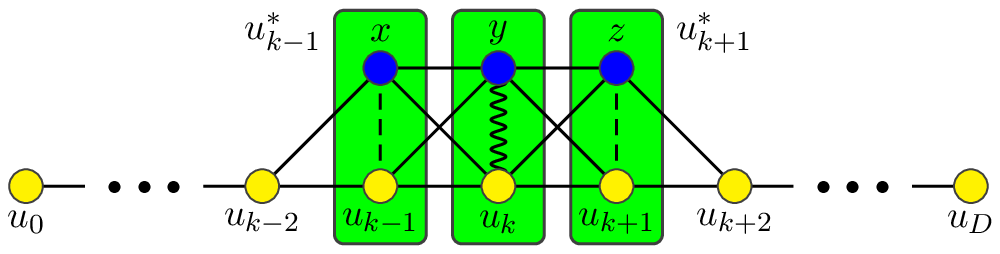}}{$\{x,y,z\}$ does not resolve $u_{k-1},u_{k+1}$ in Case 1(d).}


\textbf{Case 2.} $G^*\isomorphic P_{D+1,k}$ for some $k\in[3,D-1]$: Thus $G^*$ is path $(u_0^*,u_1^*,\dots,u_D^*)$ plus one vertex $w^*$ adjacent to $u_{k-1}^*$. As illustrated in \figref{Figure22}, suppose that every vertex of $G^*$ is of type $(1)$, except for $u_{k-2}^*$, $u_{k}^*$ and $w^*$ which are type $(1N)$, and $u_{k-1}^*$ which is of any type. In this case $n-D-1=\alpha (G^*)+1$. Consequently, it suffices to prove that $\alpha(G^*)+1$ vertices do not resolve $G$. Suppose there is a resolving set $W$ in $G$ of cardinality $\alpha(G^*)+1$. By \corref{totTwinMenysUn}, we can assume that $W$ contains the $\alpha(G^*)$ twins of $u_{k-2}$, $u_{k-1}$, $u_{k}$ and $w$ (if they exist), and another vertex of $G$. Let $x_{k-2}$, $x_{k-1}$, $x_{k}$ and $y$ respectively be twin vertices of $u_{k-2}$, $u_{k-1}$, $u_{k}$ and $w$ (if they exist). Then $x_{k-2}\not\sim u_{k-2}$,  $x_{k}\not\sim u_{k}$, and $y\not\sim w$. Thus the distance from $x_{k-2}$ (respectively $x_{k-1}$, $x_{k}$, $y$) to any vertex of $u_{k-2}$, $u_{k}$, $w$ is $2$ (respectively $1$, $2$, $2$). Hence any set of twins of vertices in $\{u_{k-2},u_{k-1},u_k,w\}$ (if they exist) does not resolve $\{u_{k-2},u_{k},w\}$. Moreover, if $i\in [0,k-1]$ then $u_i$ does not resolve $u_{k},w$; if $i\in[k-1,D]$ then $u_i$ does not resolve $u_{k-2},w$; and $w$ does not resolve $u_{k-2},u_{k}$. Therefore, $\alpha(G^*)+1$ vertices do not resolve $G$.

\Figure{Figure22}{\includegraphics{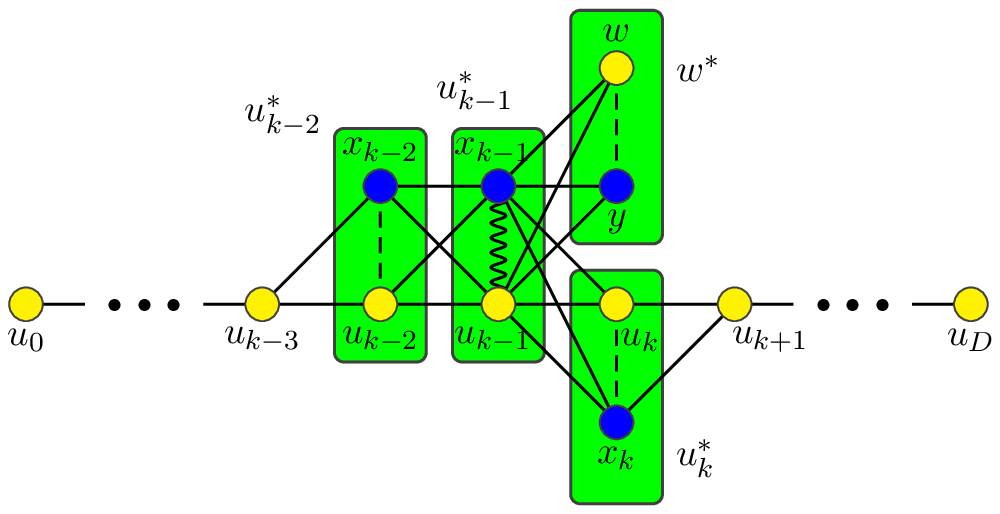}}{$\{x_{k-2},x_{k-1},x_k,y\}$ does not resolve $\{u_{k-2},u_{k},w\}$ in Case 2.}


\textbf{Case 3.} $G^*\isomorphic P'_{D+1,k}$ for some $k\in[2,D-1]$: Thus $G^*$ is path $(u_0^*,u_1^*,\dots,u_D^*)$ plus one vertex $w^*$ adjacent to $u_{k-1}^*$ and $u_{k}^*$. As illustrated in \figref{Figure23}, suppose that every vertex of $G^*$ is type $(1)$ except for $u_{k-1}^*$, $u_{k}^*$, and $w^*$ which are of type $(1K)$. In this case, $n-D-1=\alpha (G^*)+1$. Consequently, it suffices to prove that $\alpha (G^*)+1$ vertices do not resolve $G$. Suppose there is a resolving set $W$ in $G$ of cardinality $\alpha(G^*)+1$. By \corref{totTwinMenysUn}, we may assume that $W$ contains exactly the $\alpha(G^*)$ twin vertices of $u_{k-1}$, $u_{k}$ and $w$ (if they exist), and another vertex of $G$. Let $x_{k-1}$, $x_{k}$, and $y$ respectively be twins of $u_{k-1}$, $u_{k}$ and $w$ (if they exist). Hence $x_{k-1}\sim u_{k-1}$, $x_{k}\sim u_{k}$ and $y\sim w$. Consequently, $u_{k-1}$, $u_{k}$ and $w$ are at distance $1$ from $x_k$, $x_{k+1}$ and $y$. Thus any set of twins of vertices in $\{u_{k-1},u_{k},w\}$ (if they exist) does not resolve $\{u_{k-1},u_{k},w\}$. Moreover, if $i\in [0,k-1]$ then $u_i$ does not resolve $u_{k},w$; if $i\in [k,D]$ then $u_i$ does not resolve $u_{k-1},w$; and $w$ does not resolve $u_{k-1},u_{k}$. Thus $\alpha(G^*)+1$ vertices do not resolve $G$.\qed

\Figure{Figure23}{\includegraphics{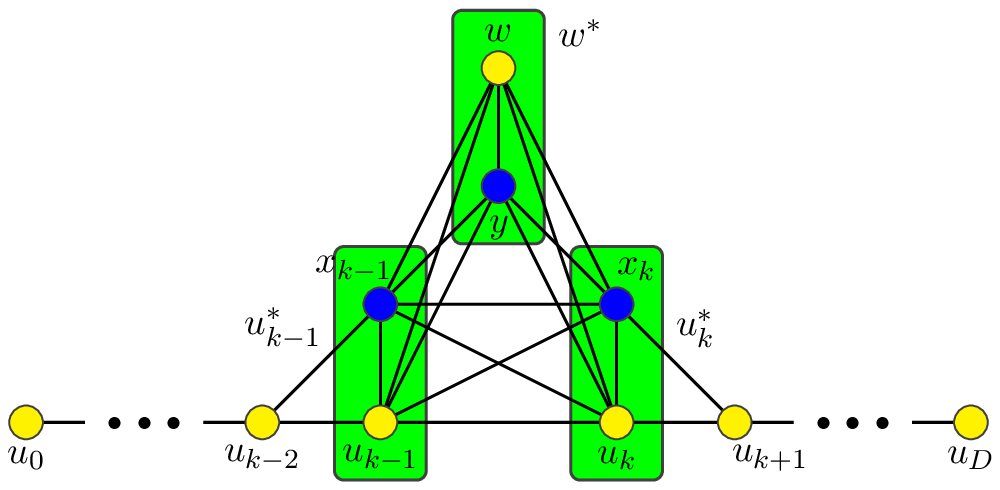}}{$\{x_{k-1},x_{k},y\}$ does not resolve $\{u_{k-1},u_{k},w\}$ in Case 3.}

\SECT{MaxOrder}{Graphs with Maximum Order}

In this section we determine the maximum order of a graph in \G.

\THM{MaxOrder}{For all integers $D\geq 2$ and $\beta\geq 1$, the maximum order of a connected graph with diameter $D$ and metric dimension $\beta$ is
\EQN{MaxOrder}{\bracket{\FLOOR{\frac{2D}{3}}+1}^\beta +
\beta\sum_{i=1}^{\CEIL{D/3}}(2i-1)^{\beta-1}\enspace.}}

First we prove the upper bound in \thmref{MaxOrder}.

\LEM{UpperBound}{For every graph $G\in\G$,
\begin{equation*}
|V(G)|\leq \bracket{\FLOOR{\frac{2D}{3}}+1}^\beta +
\beta\sum_{i=1}^{\CEIL{D/3}}(2i-1)^{\beta-1}\enspace.
\end{equation*}}

\PROOF{Let $S$ be a metric basis of $G$.
Let $k\in[0,D]$ be specified later.
For each vertex $v\in S$ and integer $i\in[0,k]$, let $N_i(v):=\{x\in V(G):\dist(v,x)=i\}$.

Consider two vertices $x,y\in N_i(v)$.
There is a path from $x$ to $v$ of length $i$, and
there is a path from $y$ to $v$ of length $i$.
Thus $\dist(x,y)\leq 2i$.
Hence for each vertex $u\in S$,
the difference between $\dist(u,x)$ and $\dist(u,y)$ is at most $2i$.
Thus the distance vector of $x$ with respect to $S$ has an $i$ in the coordinate corresponding to $v$, and in each other coordinate, there are at most $2i+1$ possible values.
Therefore $|N_i(v)|\leq (2i+1)^{\beta-1}$.

Consider a vertex $x\in V(G)$ that is not in $N_i(v)$ for all $v\in S$ and $i\in[0,k]$. Then $\dist(x,v)\geq k+1$ for all $v\in S$. Thus the distance vector of $x$ with respect to $S$ consists of $\beta$ numbers in $[k+1,D]$. Thus there are at most $(D-k)^\beta$ such vertices.  Hence
\begin{align*}
|V(G)|
&\leq(D-k)^\beta+\sum_{v\in S}\sum_{i=0}^{k}|N_i(v)|\\
&\leq (D-k)^\beta+ \beta\sum_{i=0}^{k}(2i+1)^{\beta-1}\enspace.
\end{align*}
Note that with $k=0$ we obtain the bound $|V(G)|\leq D^\beta+\beta$, independently due to \citet{KRR-DAM96} and \citet{CEJO-DAM00}. Instead we define $k:=\ceil{D/3}-1$. Then $k\in[0,D]$ and
\begin{align*}
|V(G)|
&\leq \BRACKET{D-\CEIL{\frac{D}{3}}+1}^\beta+ \beta\sum_{i=0}^{\ceil{D/3}-1}(2i+1)^{\beta-1}\\
&= \BRACKET{\FLOOR{\frac{2D}{3}}+1}^\beta+\beta\sum_{i=1}^{\ceil{D/3}}(2i-1)^{\beta-1}
\enspace.
\end{align*}
}



To prove the lower bound in \thmref{MaxOrder} we construct a graph $G\in\G$ with as many vertices as in \eqnref{MaxOrder}. The following definitions apply for the remainder of this section. Let $A:=\ceil{D/3}$ and $B:=\ceil{D/3}+\floor{D/3}$. Consider the following subsets of $\mathbb{Z}^\beta$. Let
\begin{equation*}
Q:=\{(x_1,\dots,x_\beta):A\leq x_i\leq D,i\in[1,\beta]\}\enspace.
\end{equation*}
For each $i\in[1,\beta]$ and $r\in[0,A-1]$, let
\begin{equation*}
P_{i,r}:=\{(x_1,\dots,x_{i-1},r,x_{i+1},\dots ,x_\beta):x_j\in[B-r,B+r],j\ne i \}\enspace.
\end{equation*}
Let $P_i:=\bigcup\{P_{i,r}:r\in[0,A-1]\}$ and
$P:=\bigcup\{P_i:i\in[1,\beta]\}$. Let $G$ be the graph with vertex set $V(G):=Q\cup P$, where two vertices $(x_1,\dots,x_\beta)$ and $(y_1,\dots,y_\beta)$ in $V(G)$ are adjacent if and only if $|y_i-x_i|\leq 1$ for each $i\in[1,\beta]$. \twofigref{MaxOrder}{Shape} illustrate $G$ for $\beta=2$ and $\beta=3$ respectively.

\Figure{MaxOrder}{\includegraphics{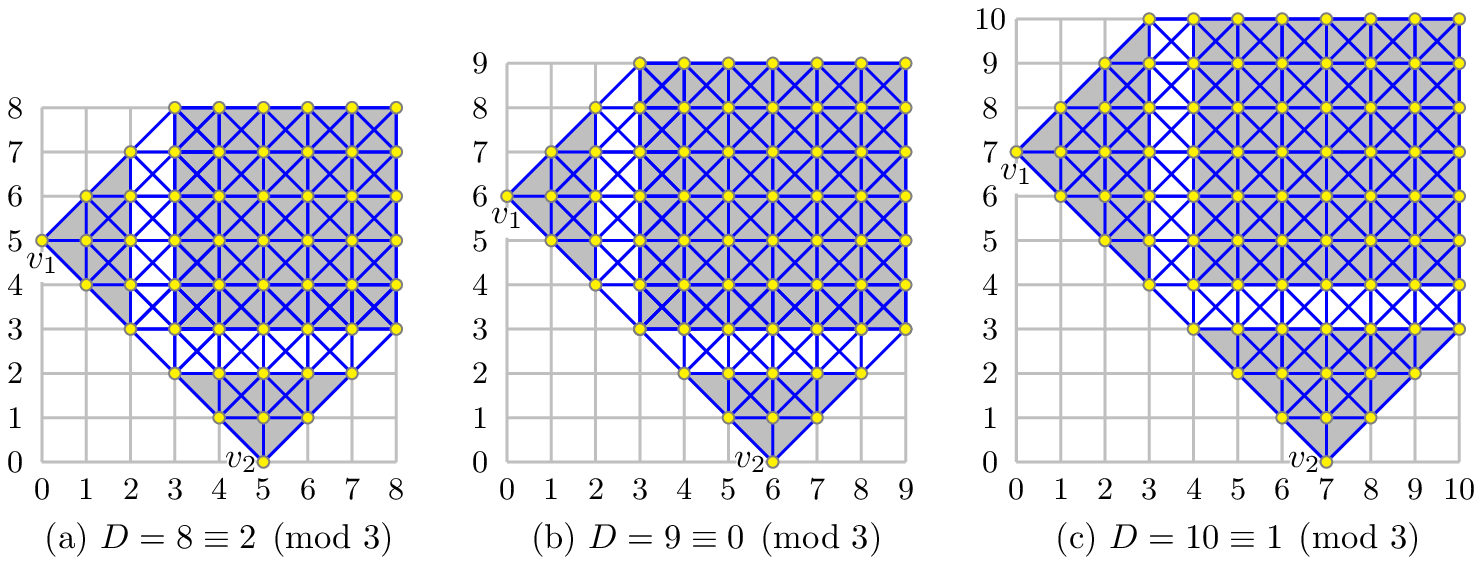}}{The graph $G$ with $\beta=2$. The shaded regions are $Q$, $P_1$, and $P_2$.}

\Figure{Shape}{\includegraphics[clip=true,bb=275 420 527 660, scale=0.5]{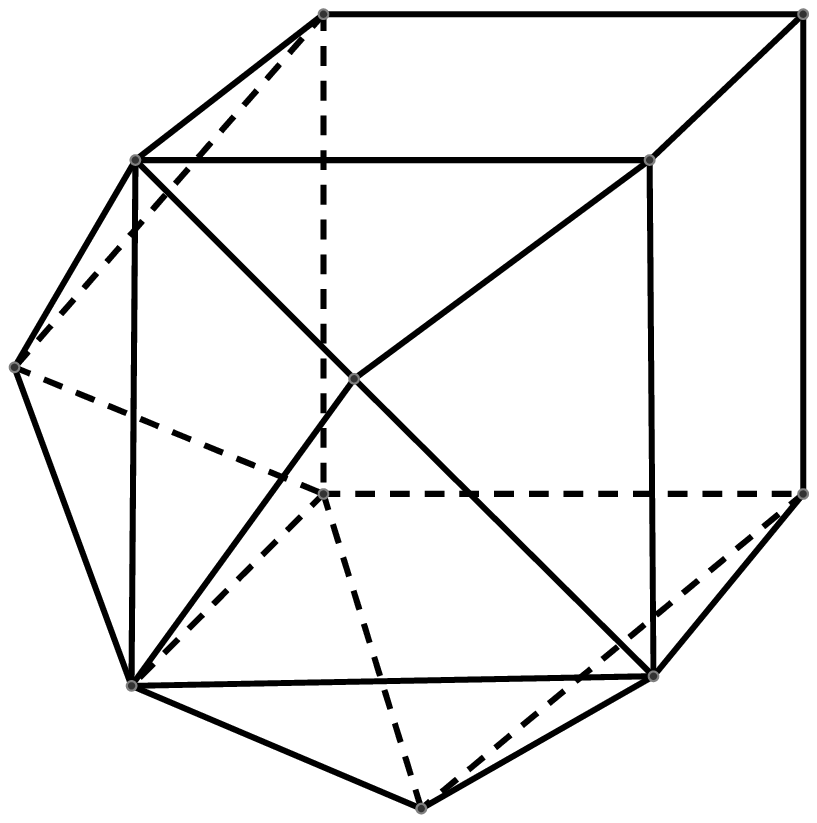}}{The convex hull of $V(G)$ with $\beta=3$.}

\LEM{OrderG}{For all positive integers $D$ and $\beta$, 
\begin{equation*}
|V(G)|=\bracket{\FLOOR{\frac{2D}{3}}+1}^\beta +
\beta\sum_{i=1}^{\CEIL{D/3}}(2i-1)^{\beta-1}\enspace.
\end{equation*}}

\PROOF{Observe that each coordinate of each vertex in $Q$ is at least $A$, and each vertex in $P$ has some coordinate less than $A$. Thus $Q\cap P=\emptyset$. Each vertex in $P_j$ ($j\ne i$) has an $i$-coordinate at least $B-r\geq B-(A-1)=\floor{D/3}+1$, and each vertex in $P_i$ has an $i$-coordinate of $r\leq A-1<\floor{D/3}+1$. Thus $P_i\cap P_j=\emptyset$ whenever $i\ne j$. Each vertex in $P_{i,r}$ has an $i$-coordinate of $r$. Thus $P_{i,r}\cap P_{i,s}=\emptyset$ whenever $r\ne s$. Thus
\begin{align*}
|V(G)|
&=|Q|+\sum_{i=1}^\beta\sum_{r=0}^{A-1}|P_{i,r}|\\
&=\BRACKET{D-\BRACKET{\CEIL{\frac{D}{3}}-1}}^\beta + \beta\sum_{r=0}^{A-1}(2r+1)^{\beta-1}\\
&=\BRACKET{\FLOOR{\frac{2D}{3}}+1}^\beta + \beta\sum_{r=1}^{A}(2r-1)^{\beta-1}
\enspace.
\end{align*}
}



We now determine the diameter of $G$. For distinct vertices $x=(x_1,\dots,x_\beta)$ and $y=(y_1,\dots,y_\beta)$ of $G$, let $z(x,y):=(z_1,\dots,z_\beta)$ where
\begin{equation*}
z_i= \begin{cases}
 x_i & \text{if }x_i=y_i, \\
 x_i+1 & \text{if }x_i<y_i, \\
 x_i-1 & \text{if }x_i>y_i\enspace.
\end{cases}
\end{equation*}

\LEM{siguiente}{$z(x,y)\in V(G)$ for all distinct vertices $x,y\in V(G)$.}

\PROOF{The following observations are an immediate consequence of the definition of $z(x,y)$, where $h,k\in \mathbb{Z}$ and $j\in[1,\beta]$:
\begin{enumerate}
\item[(i)] if $x_j,y_j\in[h,k]$ then $z_j\in[h,k]$;
\item[(ii)] if $x_j\in[h,k]$ then   $z_j\in [h-1,k+1]$; and
\item[(iii)] if $x_j\in[h,k]$ and $y_j\in[h',k']$ for some $h'>h$ and $k'<k$, then $z_j\in[h+1,k-1]$.
\end{enumerate}

We distinguish the following cases:

\begin{itemize}
\item[(1)] $x,y\in Q$: Then $x_j,y_j\in[A,D]$ for all $j$. Thus $z_j\in[a,D]$ by (i). Hence $z\in Q$.

\item[(2)] $x\in P$ and $y\in Q$: Without loss of generality, $x\in P_{1,r}$; that is, $x=(r,x_2,\dots,x_\beta)$, where $r\in[0,A-1]$ and $x_j\in[B-r,B+r]$. Since $y\in Q$, we have $y_1\ge A>r$. Thus $z=(r+1,z_2,\dots,z_\beta)$.
\begin{itemize}
\item[(2.1)] $z_1=r+1<A$: By (ii), $z_j\in[B-r-1,B+r+1]$ for every $j\neq 1$. Thus $z\in P_{1,r+1}$.
\item[(2.2)] $z_1=r+1=A$: Then $z_1\in[A,D]$. On the other hand,
if $j\neq 2$ then $y_j\geq A$, and since $x\in P_{1,r}$, we have $r=A-1$,
we have $x_j\geq B-r=\floor{D/3}+1\geq\ceil{D/3}=A$.
Thus $z_j\ge A$. Hence $z\in Q$.
\end{itemize}

\item[(3)] $x\in Q$ and $y\in P$: Without loss of generality, $y\in P_{1,r}$. That is, $y=(r,y_2,\dots,y_\beta)$, where $r\in[0,A-1]$, $y_j\in[B-r,B+r]$, and $x=(x_1,\dots,x_\beta)$, where $x_j\geq A$ for all $j$. That is, $x_1>r=y_1$, and therefore $z=(x_1-1,z_2,\dots,z_\beta)$.
\begin{itemize}
\item[(3.1)] $z_1=x_1-1\ge A$:  Since $x_j,y_j\ge A$, by (i), $z_j\ge A$ for every $j\neq 1$. That is, $z\in Q$.
\item[(3.2)] $z_1=x_1-1=A-1$: Since $r\leq A-1$, we have $y_j\in[B-r,B+r]\subseteq[B-(A-1),B+A-1]$ for all $j\not\in\{1,2\}$. Now $B-A=\floor{D/3}\leq\ceil{D/3}=A$ and $D\leq\floor{D/3}+2\ceil{D/3}=A+B$. Thus 
$x_j\in[B-A,B+A]$. By (iii), $z_j\in [B-(A-1),B+A-1]$. That is, $z\in P_{1,A-1}$.
\end{itemize}

\item[(4)] $x,y\in P_h$: Without loss of generality, $x,y\in P_1$. Thus $x=(r,x_2,\dots,x_\beta)$ for some $r\in[0,A-1]$ with $x_j\in[B-r,B+r]$ for all $j\neq1$, and $y=(s,y_2,\dots,y_\beta)$ for some $s\in[0,A-1]$ with $y_j\in[B-s,B+s]$ for all $j\neq1$.
\begin{itemize}
\item[(4.1)] $r=s$: Then $z=(r,z_2,\dots ,z_\beta )$. 
By (i), $z_j\in[B-r,B+r]$ for all $j\neq1$. Thus $z\in P_{1,r}$.
\item[(4.2)] $r<s$: Then $z=(r+1,z_2,\dots,z_\beta)$. By (ii), $z_j\in[B-(r+1),B+r+1]$ for all $j\neq1$. Thus $z\in P_{1,r+1}$. 
\item[(4.3)] $r>s$: Then $z=(r-1,z_2,\dots,z_\beta)$. By (iii),   $z_j\in[B-(r-1),B+r-1]$. Thus $z\in P_{1,r-1}$.
\end{itemize}

\item[(5)] $x\in P_h$, $y\in P_k$ and $h\neq k$: Without loss of generality, $x\in P_1$ and $y\in P_2$. Thus $x=(r,x_2,\dots,x_\beta)$ for some $r\in[0,A-1]$ with $x_j\in[B-r,B+r]$ for all $j\neq1$, and $y=(y_1,s,y_3,\dots,y_\beta)$ for some $s\in[0,A-1]$ with $y_j\in[B-s,B+s]$ for all $j\neq2$. Hence $r<A\leq y_1$ and $s<A\leq x_2$, implying $z=(r+1,x_2-1,z_3,\dots,z_\beta)$.
\begin{itemize}
\item[(5.1)] $z_1=r+1<A$: Now $x_j\in[B-r,B+r]$ for $j\neq 1$. Thus $z_j\in[B-r-1,B+r+1]$ by (ii). Thus $z\in P_{1,r+1}$.
\item[(5.2)] $z_1=r+1=A$: Consider the following subcases:
\begin{itemize}
\item[(5.2.1)] $z_2=x_2-1\geq A$: By hypotheses, $z_1,z_2\geq A$. For
$j\not\in\{1,2\}$, since $x_j,y_j\ge A$, (i) implies that $z_j\geq A$. Thus $z\in Q$.
\item[(5.2.2)] $z_2=x_2-1=A-1$: In this case $x=(A-1,A,x_3,\dots,x_\beta)$,
$z=(A,A-1,z_3,\dots,z_\beta)$, and $s\leq A-1=r$. Since $x\in P_{1,A-1}$, we have $z_1=x_2=A\in[B-(A-1),B+A-1]$. For $j\not\in\{1,2\}$, since $x_j\in[B-(A-1),B+A-1]$ and $y_j\in[B-s,B+s]$, where $s\leq r=A-1$, (i) implies that $z_j\in[B-(A-1),B+A-1]$. That is, $z\in P_{2,A-1}$.
\end{itemize}
\end{itemize}
\end{itemize}
}

\LEM{WeNeedToProve}{For all vertices $x=(x_1,\dots,x_\beta)$ and $y=(y_1,\dots,y_\beta)$ of $G$,
\eqn{\dist(x,y)=\max\{|y_i-x_i|:i\in[1,\beta]\}\leq D.}}

\PROOF{ 
For each $i\in[1,\beta]$, $\dist(x,y)\geq|x_i-y_i|$ since on every $xy$-path $P$, the $i$-coordinates of each pair of adjacent vertices in $P$ differ by at most $1$. This proves the lower bound $\dist(x,y)\geq\max_i|y_i-x_i|$. 

Now we prove the upper bound $\dist(x,y)\leq\max_i|y_i-x_i|$ by induction. If $\max_i|y_i-x_i|=1$ then $x$ and $y$ are adjacent, and thus $\dist(x,y)=1$. Otherwise, let $z:=z(x,y)$. By the definition of $z(x,y)$, for all $i\in[1,\beta]$ we have $|y_i-z_i|=|y_i-x_i|-1$ unless $x_i=y_i$. Thus $\max_i|y_i-z_i|=\max_i|y_i-x_i|-1$. By induction, 
$\dist(z,y)\leq\max_i|y_i-z_i|=\max_i|y_i-x_i|-1$. By \lemref{siguiente}, $z$ is a vertex of $G$, and by construction, $x$ and $z$ are adjacent. Thus $\dist(x,y)\leq \dist(z,y)+1=\max_i|y_i-x_i|$, as desired. 
}

\lemref{WeNeedToProve} implies that $G$ has diameter $D$. Let $S:=\{v_1,\dots,v_\beta \}$, where 
\eqn{v_i=(\underbrace{t,\dots,t}_{i-1},0,t,\dots t)\enspace.} 
Observe that each $v_i\in P_i$. We now prove that $S$ is a metric basis of $G$.

\LEM{distxvi}{$\dist(x,v_i)=x_i$ for every vertex $x=(x_1,\dots ,x_\beta)$ of $G$ and for each $v_i\in S$.}

\PROOF{Let $v_{i,j}$ be the $j$-th coordinate of $v_i$; that is, $v_{i,i}=0$ and $v_{i,j}=t$ for $i\ne j$. Then
$\dist(x,v_i)
=\max \{ |v_{i,j}-x_j|:1\le j \le \beta\} =\max \{x_i,\max\{|t-x_j|:1\le j \le \beta,j\ne i\}\}$.
We claim that $|t-x_j|\leq x_i$ for each $j\ne i$, implying $\dist(x,v_i)=x_i$, as desired.

First suppose that $x\in Q$. Then $s\leq x_j\leq D$. Thus
$|B-x_j|
\leq\newline\max\{B-A,D-t\}
=\max\{\floor{D/3},D-t\}
\leq\max\{\floor{D/3},\ceil{D/3}\}
=\ceil{D/3}
\leq x_i$.

Now suppose that $x\in P_{k,r}$ for some $k\ne i$ and for some $r$. Then $x_i\geq
t-r\geq t-(\ceil{D/3}-1)=\floor{D/3}+1\geq\ceil{D/3}$.
Now $|t-x_j|\leq r\leq \ceil{D/3}-1$. Thus $|t-x_j|\leq x_i$.

Finally suppose that $x\in P_{i,r}$ for some $r$. Then $|t-x_j|\leq r=x_i$.}

\lemref{distxvi} implies that the metric coordinates of a vertex $x\in V(G)$ with respect to $S$ are its coordinates as elements of $\mathbb{Z}^\beta$. Therefore $S$ resolves $G$. Thus $G$ has metric dimension at most $|S|=\beta$.

If the metric dimension of $G$ was less than $\beta$, then by \lemref{UpperBound},
\begin{equation*}
\bracket{\FLOOR{\frac{2D}{3}}+1}^{\beta} +
\beta\sum_{i=1}^{\CEIL{D/3}}(2i-1)^{\beta-1}
= |V(G)|\leq
\bracket{\FLOOR{\frac{2D}{3}}+1}^{\beta-1} +
(\beta-1)\sum_{i=1}^{\CEIL{D/3}}(2i-1)^{\beta-2}\enspace,
\end{equation*}
which is a contradiction. Thus $G$ has metric dimension $\beta$, and $G\in\G$. This completes the proof of \thmref{MaxOrder}.


\begin{thebibliography}{33}
\providecommand{\natexlab}[1]{#1}
\providecommand{\url}[1]{\texttt{#1}}
\providecommand{\urlprefix}{}
\expandafter\ifx\csname urlstyle\endcsname\relax
  \providecommand{\doi}[1]{doi:\discretionary{}{}{}#1}\else
  \providecommand{\doi}{doi:\discretionary{}{}{}\begingroup
  \urlstyle{rm}\Url}\fi

\bibitem[{Beerliova et~al.(2006)Beerliova, Eberhard, Erlebach, Hall, Hoffmann,
  Mihal'\'ak, and Ram}]{BEEHHMS-IEEE06}
\textsc{Zuzana Beerliova, Felix Eberhard, Thomas Erlebach, Alexander Hall,
  Michael Hoffmann, Mat\'u\v{s} Mihal'\'ak, and L.~Shankar Ram}.
\newblock Network discovery and verification.
\newblock \emph{IEEE J. on Selected Areas in Communications},
  24(12):2168--2181, 2006.

\bibitem[{Brigham et~al.(2003)Brigham, Chartrand, Dutton, and
  Zhang}]{BCDZ-MB03}
\textsc{Robert~C. Brigham, Gary Chartrand, Ronald~D. Dutton, and Ping Zhang}.
\newblock Resolving domination in graphs.
\newblock \emph{Math. Bohem.}, 128(1):25--36, 2003.

\bibitem[{C{\'a}ceres et~al.(pear)C{\'a}ceres, Hernando, Mora, Pelayo, Puertas,
  Seara, and Wood}]{Caceres-etal}
\textsc{Jos{\'e} C{\'a}ceres, Carmen Hernando, Merc{\`e} Mora, Ignacio~M.
  Pelayo, Mar{\'\i}a~L. Puertas, Carlos Seara, and David~R. Wood}.
\newblock On the metric dimension of cartesian products of graphs.
\newblock \emph{SIAM J. Discrete Math.}, to appear.
\newblock \urlprefix\url{arXiv.org/math/0507527}.

\bibitem[{Chappell et~al.(2003)Chappell, Gimbel, and Hartman}]{CGH}
\textsc{Glenn~G. Chappell, John Gimbel, and Chris Hartman}.
\newblock Bounds on the metric and partition dimensions of a graph, 2003.
\newblock
  \urlprefix\url{http://www.cs.uaf.edu/{\textasciitilde}chappell/papers/metric%
/}.

\bibitem[{Chartrand et~al.(2000{\natexlab{a}})Chartrand, Eroh, Johnson, and
  Oellermann}]{CEJO-DAM00}
\textsc{Gary Chartrand, Linda Eroh, Mark~A. Johnson, and Ortrud~R. Oellermann}.
\newblock Resolvability in graphs and the metric dimension of a graph.
\newblock \emph{Discrete Appl. Math.}, 105(1-3):99--113, 2000{\natexlab{a}}.

\bibitem[{Chartrand et~al.(2000{\natexlab{b}})Chartrand, Poisson, and
  Zhang}]{CPZ-CMA00}
\textsc{Gary Chartrand, Christopher Poisson, and Ping Zhang}.
\newblock Resolvability and the upper dimension of graphs.
\newblock \emph{Comput. Math. Appl.}, 39(12):19--28, 2000{\natexlab{b}}.

\bibitem[{Chartrand and Zhang(2003)}]{CZ-CN03}
\textsc{Gary Chartrand and Ping Zhang}.
\newblock The theory and applications of resolvability in graphs. {A} survey.
\newblock In \emph{Proc. 34th Southeastern International Conf. on
  Combinatorics, Graph Theory and Computing}, vol. 160 of \emph{Congr. Numer.},
  pp. 47--68. 2003.

\bibitem[{Chv{\'a}tal(1983)}]{Chvatal-Comb83}
\textsc{Va{\v{s}}ek Chv{\'a}tal}.
\newblock Mastermind.
\newblock \emph{Combinatorica}, 3(3-4):325--329, 1983.

\bibitem[{Currie and Oellermann(2001)}]{CO01}
\textsc{James Currie and Ortrud~R. Oellermann}.
\newblock The metric dimension and metric independence of a graph.
\newblock \emph{J. Combin. Math. Combin. Comput.}, 39:157--167, 2001.

\bibitem[{Erd{\H{o}}s and R{\'e}nyi(1963)}]{ER63}
\textsc{Paul Erd{\H{o}}s and Alfr{\'e}d R{\'e}nyi}.
\newblock On two problems of information theory.
\newblock \emph{Magyar Tud. Akad. Mat. Kutat\'o Int. K\"ozl.}, 8:229--243,
  1963.

\bibitem[{Goddard(2003)}]{Goddard-JCMCC03}
\textsc{Wayne Goddard}.
\newblock Static mastermind.
\newblock \emph{J. Combin. Math. Combin. Comput.}, 47:225--236, 2003.

\bibitem[{Goddard(2004)}]{Goddard-JCMCC04}
\textsc{Wayne Goddard}.
\newblock Mastermind revisited.
\newblock \emph{J. Combin. Math. Combin. Comput.}, 51:215--220, 2004.

\bibitem[{Greenwell(2000)}]{Greenwell-JRM99}
\textsc{Don~L. Greenwell}.
\newblock Mastermind.
\newblock \emph{J. Recr. Math.}, 30:191--192, 1999-2000.

\bibitem[{Guy and Nowakowski(1995)}]{GN-AMM95}
\textsc{Richard~K. Guy and Richard~J. Nowakowski}.
\newblock Coin-weighing problems.
\newblock \emph{Amer. Math. Monthly}, 102(2):164, 1995.

\bibitem[{Harary and Melter(1976)}]{HM-AC76}
\textsc{Frank Harary and Robert~A. Melter}.
\newblock On the metric dimension of a graph.
\newblock \emph{Ars Combinatoria}, 2:191--195, 1976.

\bibitem[{Kabatianski et~al.(2000)Kabatianski, Lebedev, and Thorpe}]{KLT00}
\textsc{Grigori Kabatianski, V.~S. Lebedev, and J.~Thorpe}.
\newblock The {M}astermind game and the rigidity of {H}amming spaces.
\newblock In \emph{Proc. IEEE International Symposium on Information Theory
  \textup{(ISIT '00)}}, p. 375. IEEE, 2000.

\bibitem[{Khuller et~al.(1996)Khuller, Raghavachari, and Rosenfeld}]{KRR-DAM96}
\textsc{Samir Khuller, Balaji Raghavachari, and Azriel Rosenfeld}.
\newblock Landmarks in graphs.
\newblock \emph{Discrete Appl. Math.}, 70(3):217--229, 1996.

\bibitem[{Lindstr{\"o}m(1964)}]{Lindstrom64}
\textsc{Bernt Lindstr{\"o}m}.
\newblock On a combinatory detection problem. {I}.
\newblock \emph{Magyar Tud. Akad. Mat. Kutat\'o Int. K\"ozl.}, 9:195--207,
  1964.

\bibitem[{Manuel et~al.(2006)Manuel, Rajan, Rajasingh, and
  Chris~Monica}]{MRRC-JDMSC06}
\textsc{Paul Manuel, Bharati Rajan, Indra Rajasingh, and M.~Chris~Monica}.
\newblock Landmarks in torus networks.
\newblock \emph{J. Discrete Math. Sci. Cryptogr.}, 9(2):263--271, 2006.

\bibitem[{Peters-Fransen and Oellermann(2006)}]{PO-UM06}
\textsc{Joel Peters-Fransen and Ortrud~R. Oellermann}.
\newblock The metric dimension of the cartesian product of graphs.
\newblock \emph{Utilitas Math.}, 69:33--41, 2006.

\bibitem[{Poisson and Zhang(2002)}]{PZ02}
\textsc{Christopher Poisson and Ping Zhang}.
\newblock The metric dimension of unicyclic graphs.
\newblock \emph{J. Combin. Math. Combin. Comput.}, 40:17--32, 2002.

\bibitem[{Saenpholphat and Zhang(2003)}]{SZ-CMJ03}
\textsc{Varaporn Saenpholphat and Ping Zhang}.
\newblock Connected resolvability of graphs.
\newblock \emph{Czechoslovak Math. J.}, 53(128)(4):827--840, 2003.

\bibitem[{Saenpholphat and Zhang(2004{\natexlab{a}})}]{SZ-IJMMS04}
\textsc{Varaporn Saenpholphat and Ping Zhang}.
\newblock Conditional resolvability in graphs: a survey.
\newblock \emph{Int. J. Math. Math. Sci.}, (37-40):1997--2017,
  2004{\natexlab{a}}.

\bibitem[{Saenpholphat and Zhang(2004{\natexlab{b}})}]{SZ-CN04}
\textsc{Varaporn Saenpholphat and Ping Zhang}.
\newblock Detour resolvability of graphs.
\newblock In \emph{Proc. of 35th Southeastern International Conf. on
  Combinatorics, Graph Theory and Computing}, vol. 169 of \emph{Congr. Numer.}, pp. 3--21. Utilitas Math., 2004{\natexlab{b}}.

\bibitem[{Saenpholphat and Zhang(2004{\natexlab{c}})}]{SZ-CMJ04}
\textsc{Varaporn Saenpholphat and Ping Zhang}.
\newblock On connected resolving decompositions in graphs.
\newblock \emph{Czechoslovak Math. J.}, 54(129)(3):681--696,
  2004{\natexlab{c}}.

\bibitem[{Seb{\H{o}} and Tannier(2004)}]{ST-MOR04}
\textsc{Andr{\'a}s Seb{\H{o}} and Eric Tannier}.
\newblock On metric generators of graphs.
\newblock \emph{Math. Oper. Res.}, 29(2):383--393, 2004.

\bibitem[{Shanmukha et~al.(2002)Shanmukha, Sooryanarayana, and
  Harinath}]{SSH02}
\textsc{B.~Shanmukha, B.~Sooryanarayana, and K.~S. Harinath}.
\newblock Metric dimension of wheels.
\newblock \emph{Far East J. Appl. Math.}, 8(3):217--229, 2002.

\bibitem[{Slater(1975)}]{Slater75}
\textsc{Peter~J. Slater}.
\newblock Leaves of trees.
\newblock In \emph{Proc. 6th Southeastern Conf. on Combinatorics, Graph Theory,
  and Computing}, vol.~14 of \emph{Congr. Numer.}, pp. 549--559. Utilitas
  Math., 1975.

\bibitem[{Slater(1988)}]{Slater-JMPS88}
\textsc{Peter~J. Slater}.
\newblock Dominating and reference sets in a graph.
\newblock \emph{J. Math. Phys. Sci.}, 22(4):445--455, 1988.

\bibitem[{S{\"o}derberg and Shapiro(1963)}]{SS-AMM63}
\textsc{Staffan S{\"o}derberg and Harold~S. Shapiro}.
\newblock A combinatory detection problem.
\newblock \emph{Amer. Math. Monthly}, 70:1066, 1963.

\bibitem[{Sooryanarayana(1998)}]{Sooryanarayana98}
\textsc{B.~Sooryanarayana}.
\newblock On the metric dimension of a graph.
\newblock \emph{Indian J. Pure Appl. Math.}, 29(4):413--415, 1998.

\bibitem[{Sooryanarayana and Shanmukha(2001)}]{SS01}
\textsc{B.~Sooryanarayana and B.~Shanmukha}.
\newblock A note on metric dimension.
\newblock \emph{Far East J. Appl. Math.}, 5(3):331--339, 2001.

\bibitem[{Yushmanov(1987)}]{Yushmanov87}
\textsc{S.~V. Yushmanov}.
\newblock Estimates for the metric dimension of a graph in terms of the
  diameters and the number of vertices.
\newblock \emph{Vestnik Moskov. Univ. Ser. I Mat. Mekh.}, 103:68--70, 1987.

\end{thebibliography}

\def\soft#1{\leavevmode\setbox0=\hbox{h}\dimen7=\ht0\advance \dimen7
  by-1ex\relax\if t#1\relax\rlap{\raise.6\dimen7
  \hbox{\kern.3ex\char'47}}#1\relax\else\if T#1\relax
  \rlap{\raise.5\dimen7\hbox{\kern1.3ex\char'47}}#1\relax \else\if
  d#1\relax\rlap{\raise.5\dimen7\hbox{\kern.9ex \char'47}}#1\relax\else\if
  D#1\relax\rlap{\raise.5\dimen7 \hbox{\kern1.4ex\char'47}}#1\relax\else\if
  l#1\relax \rlap{\raise.5\dimen7\hbox{\kern.4ex\char'47}}#1\relax \else\if
  L#1\relax\rlap{\raise.5\dimen7\hbox{\kern.7ex
  \char'47}}#1\relax\else\message{accent \string\soft \space #1 not
  defined!}#1\relax\fi\fi\fi\fi\fi\fi} \def\cprime{$'$}
  \def\Dbar{\leavevmode\lower.6ex\hbox to 0pt{\hskip-.23ex \accent"16\hss}D}

\end{document}